\documentclass[11pt,leqno]{amsart}
\usepackage{pstricks}
\usepackage{amsmath}
\usepackage{amsfonts}
\usepackage{amssymb}
\usepackage{amsthm}
\usepackage{tikz}
\usepackage{tcolorbox}
\usepackage{pgf,tikz,pgfplots}
\usepackage{mathrsfs}
\usetikzlibrary{arrows}
\usetikzlibrary{arrows}
\usepackage{float} 

\usepackage{marginnote}

\title{The Furstenberg set and its random version}

\author{Aihua Fan}
\address{School of mathematics and Statistics, Central China Normal University,  430079 Wuhan, China \& LAMFA, UMR 7352 CNRS, University of Picardie, 33 rue Saint Leu,80039 Amiens, France}
\email{ai-hua.fan@u-picardie.fr}

\author{Herv\'e Queff\'elec}
\address{CNRS, Laboratoire Paul Painlev\'e, UMR 8524 \& Labex CEMPI (ANR-LABX-0007-01),
Universit\'e de Lille                          
Cit\'e Scientifique, B\^t. M2                                  
59655 Villeneuve d'Ascq Cedex,
FRANCE }
\email{herve.queffelec@univ-lille.fr}

\author{Martine Queff\'elec}
\address{CNRS, Laboratoire Paul Painlev\'e, UMR 8524, \& Labex CEMPI (ANR-LABX-0007-01),
Universit\'e de Lille                          
Cit\'e Scientifique, B\^t. M2                                  
59655 Villeneuve d'Ascq Cedex,
FRANCE}
\email{martine.queffelec@univ-lille.fr}

\begin{document}

\begin{abstract}
We study some number-theoretic, ergodic and harmonic analysis properties of the Furstenberg set of integers $S=\{2^{m}3^{n}\}$ and compare them to those of its random analogue $T$. In this half-expository work, we show for example that $S$ is  ``Khinchin distributed", is far from being Hartman uniformly distributed while $T$ is, also that $S$ is a $\Lambda(p)$-set for all $2<p<\infty$ and that $T$ is a $p$-Rider set for all $p$ such that $4/3<p<2$. Measure-theoretic and probabilistic techniques, notably martingales, play an important role in this work.
\end{abstract}  

\subjclass[2010]{37A44, 43A46, 60G46.}
\keywords{Furstenberg set, Sidon set, Khintchin class, Uniform distribution, Martingale.}

\maketitle
\newenvironment{demo}{\medbreak\begingroup\noindent{\bf Proof. }}{\hfill$\diamondsuit$\endgroup\goodbreak\medskip}

\newtheorem{defi}{Definition}[section]
\newtheorem{theo}{Theorem}[section]
\newtheorem{prop}[theo]{Proposition}
\newtheorem{cor}[theo]{Corollary}
\newtheorem{lem}[theo]{Lemma}
\newtheorem{rem}{Remark}
\newtheorem{ex}{Example}
\newtheorem{pb}{Problem}
\newcommand {\N}{{\mathbb{N}}} 
\newcommand {\Z}{{\mathbb{Z}}}
\newcommand {\E}{{\mathbb{E}}}
\newcommand {\R}{{\mathbb{R}}} 
\newcommand {\Q}{{\mathbb{Q}}} 
\newcommand {\C}{{\mathbb{C}}} 
\newcommand {\D}{{\cal D}} 
\newcommand {\T}{{\mathbb{T}}}

\newcommand {\F}{{\mathcal{F}}}  

\newcommand {\G}{{\mathcal{G}}}  
\renewcommand {\S}{{\mathcal{S}}}  
\newcommand {\A}{{\cal{A}}}  
\newcommand {\J}{{\cal{J}}}  
\newcommand {\B}{{\bf{BAD}}}  
\renewcommand {\d}{{\cal{D}}} 
\renewcommand {\l}{{\cal{L}}} 
\renewcommand {\r}{{\cal{R}}} 
\renewcommand {\b}{{\cal{B}}} 
\renewcommand {\k}{{\cal{K}}} 
\newcommand {\q}{{\cal{Q}}} 
\renewcommand {\P}{{\mathbb{P}}} 
\newcommand {\p}{{\cal{P}}} 
\newcommand {\noi}{\noindent}
\newcommand {\eps}{\varepsilon}
\tableofcontents

\section{Introduction}
The  goal of this paper\footnote{This paper will appear in the journal L'Enseignement Math\'ematique. Theorem 2.10 discussing the exact order of decay  of Fourier coefficients is added in this ArXiv version after the acceptance of the paper by the journal. Section 6  included here (but not in the version of  L'Enseignement Math\'ematique) would facilitate the reading of the part on Bohr topology.} , which is mainly a survey but also contains several original results,  is the study of two different sets of integers: the Furstenberg set $S$ and its randomized version denoted $T$.

 Since various notions (from harmonic analysis, probability, ergodic theory,  number theory and diophantine approximation, fractal geometry) are involved, we begin with a rather long introduction in which the necessary notions are described in a sketchy way. But we strive to add precise bibliographical references, or references to 
 the forthcoming  sections, where those notions are presented in more  detail.
 Complements can be found on harmonic analysis, probability and Banach spaces  in \cite{LQ}, on ergodic theory in \cite{Pe} and on distribution of sequences in \cite{KN}.
 \medskip
 
 {\bf Notations.} Here are summarized  the notations we adopt. $\N$ (resp.~$\N_0$) denotes the set of positive (resp.~non negative) integers, $\Z$  the set of  integers, $\{x\}$ the fractional part of the real number $x$. The cardinality of a finite set $A$ is denoted either by $|A|$, or by  $\#A$. Also for estimations, we make use of the classical Hardy's notation $\ll$ (bounded up to a positive multiplicative constant).
We will mainly (but not exclusively) consider  the group $G=\T=\R/\Z$ with $m$ as its Haar measure and we will denote $\|x\|= \inf_{n\in \Z} |x-n|$ for $x\in \T$. 
We put $L^p(\T)$ for $L^p(\T,m)$. For every $f\in L^{1}(\T)$ (integrable and 1-periodic), more generally  every finite measure $\mu$ on $\T$ and every $ n\in \Z$, we set: 
$$\int_\T f\ dm:=\int_0^1f(x)\ dx,\ \widehat{f}(n)=\int_{\T}f \overline{e}_ndm,\  \widehat{\mu}(n)=\int_{\T}\overline{e}_nd\mu$$

where  $e_n(x)=e(nx)$ with $e(x) =\exp(2i\pi x)$.
$\mathbb{E}$ (resp.~$V$) denotes the expectation (resp.~the variance) in a probability space. 
$\mathbb{E}(X|{\mathcal{C}})$ denotes the conditional expectation of an integrable random variable $X$ given the sub-$\sigma$-field $\mathcal{C}$.  
\medskip

\subsection{The set $S$}
Our main object is the Furstenberg set $S$ of positive integers defined as
$$S=\{2^j3^k\ :\  j,k\in \N_0\}.$$
This is the multiplicative semigroup of $\N$,  generated by $2$ and $3$.  More generally, for  a  finite set of  coprime numbers $\{q_1,\dots, q_s\}$, we denote by  $S(q_1,\dots, q_s)$ the multiplicative semi-group generated by  $q_1,\dots, q_s$.
If $ 2=p_1<p_2<\cdots< p_s<\cdots$ is the ordered sequence of the prime numbers,  $S(p_1,\dots, p_s)$ is nothing but $D_s:=\{n : P^{+}(n)\leq p_s \}$
where  $P^{+}(n)$  denotes 
the greatest prime divisor of $n\geq 2$. This set $D_s$ is   called the set of \textit{$s$-friable integers} and has been  intensively studied in  analytic  
number theory (see \cite{Ma, Ti1,Ti2, Ph, BDH, GV, FT, ABT}).
 
When $s=1$, $S(q_1)$ is the so-obtained Hadamard lacunary set $\{q_{1}^n, n\ge1\}$, which enjoys various properties, as well in harmonic analysis as in ergodic theory, that we shall discuss below. 
It is quite natural to explore the same properties for less sparse sets, e.g. $S(q_1,\dots, q_s)$ as soon as  $s\ge2$.

\smallskip
We particularly focus on the set $S=S(2, 3)$ because of its  closeness to a famous conjecture of Furstenberg, the $\times2\times3$--conjecture,
which asserts that {\it every continuous probability measure on the circle $\T$ which is both $2$-invariant  and $3$--invariant  must be the Haar measure $m$ of $\T$} (see Subsection 2.2.1). We make no progress on this conjecture in  our  
paper, but this problem acted as a motivation for revisiting the set $S$, which is clearly invariant under the multiplication by $2$ and $3$. By rearrangement the  Furstenberg set  gives rise to an increasing sequence of integers, called the ``Furstenberg sequence" (also denoted $S$, to ease  notation):
$$S=\{1=s_1<s_2<\dots<s_n<\dots\}.
$$

From different  viewpoints (harmonic analysis, arithmetics, and dynamical system), the set $S$ or the sequence $(s_n)$ is quite original and mysterious. Its  rate of growth is intermediate between the growth of a polynomial sequence (e.g. the sequence $(n^2)$ of squares) and that of a power sequence (e.g. the sequence $(2^n)$) and this rate 
is not so easy to handle with.  Ramanujan gave an amazing statement, of which Hardy had to give a full proof (\cite{Ha}).
 
The  sumset $\Sigma=\{2^j+3^k\}$ was well understood from the uniform distribution and  harmonic analysis viewpoints \cite{Bl}, for the reason that it appears  as a projection of a nice product set.  
 But the multiplicative structure of $S$ prevents us from using the same  projection technique.  

\medskip

We shall study  three types of properties of $S$:
\medskip

\noindent $\bullet$  {\bf Rate of increase and lacunarity}. 
We say that a sequence of integers $E=(u_n)$ is \emph{weakly lacunary} if $u_{n+1}-u_n\to\infty$, or  {\it Hadamard lacunary}  of ratio $\rho$ if  $\rho:=\liminf u_{n+1}/u_n>1$ (in the sequel, we shall say ``Hadamard set''). 
The set $S$ is not a  Hadamard set because $s_{n+1}/s_n\to 1$, which was exploited  by Furstenberg in \cite{Fu1}.
We shall look for asymptotic behaviours of $s_n, s_{n+1}-s_n, \frac{s_{n+1}}{s_n}-1$. Notice that estimates of such quantities for a Hadamard set are obvious. 
 For the set $S$, it is easy to prove that $I_{S}(N):=|S\cap [1,N]|\approx (\log N)^2$ whence
 $s_n\approx \exp(\sqrt{n})$. 
Going further needs sharp estimates and more work (Subsection 2.1.).
\bigskip

\noindent $\bullet$ {\bf Distribution issues}.

 Consider a set of integers
$ E=\{\lambda_1<\lambda_2<\cdots<\lambda_n<\cdots\}\subset\Z.$
 As is well-known (Hermann Weyl's theorem \cite{KN}), for almost every $x\in (0,1)$ the sequence $(\lambda_{n}x)$ is {\it equidistributed modulo $1$}, meaning that, for any $h\in \Z\setminus\{0\}$, we have 
 \begin{equation}\label{herweyl}   \frac{1}{N}\sum_{n=1}^N e(h\lambda_{n} x)\to 0.
 \end{equation}

 We also say more briefly: almost every $x$ is {\it $E$-normal}. 
 The exceptional set of those $x$ such that $(\lambda_{n}x)$ is not equidistributed mod 1 is Haar-negligible. 
 It is  interesting to study its properties.
 This exceptional set can be reduced to the singleton $\{0\}$. In other words, in some extreme cases, it can happen that
 \begin{equation}\label{hart}  \frac{1}{N}\sum_{n=1}^N e(\lambda_{n} x)\to 0,\ \hbox{for all}\ x\neq 0
 \end{equation}  and then we say that the set $E$ is {\it Hartman uniformly distributed}. 
 
We introduce two measuring tools.
 Hausdorff dimension is usually used to describe the size of  sets in a metric space. 
 Consider a subset $X$ of $\R^d$. For $s>0$ and $\varepsilon>0$, 
 we define 
 $$H_{\epsilon}^s (X)= \inf \left\{\sum_{i} |U_i|^s:  X \subset \bigcup_i U_i, |U_i| <\epsilon\right\},$$
 where $\{U_i\}$ is a countable family of subsets,   called an $\varepsilon$-cover of $X$, and $|U_i|$ denotes the diameter of $U_i$.
 The {\em $s$-dimensional Hausdorff measure} is then defined by 
 $H^{s}(X)=\sup_{\varepsilon>0}H^s_\varepsilon(X)$.
There is a critical exponent $\rho$, called {\em Hausdorff dimension} of $X$ and denoted by $\dim_H(X)$, such that  $H^{s}(X)=+\infty$  for $s<\rho$ and $H^{s}(X)=0$ for $s>\rho$. 

Observe that $0\le {\rm dim}_{H}(X)\le d$ and that $ {\rm dim}_{H}(X)= d$ as soon as $X$ has positive Lebesgue measure.  For example, if $X$ is the Cantor middle-third set,  $\hbox{\ dim}_{H} (X)=\log 2/\log 3$. 
See \cite{Fal}, chapter 2 for more details about the Hausdorff dimension. 

Next,  recall that a probability measure $\mu$ on $\mathbb{T}$ is called a  {\it Rajchman measure} if $\widehat\mu (n)\to 0$ at infinity.

 The support of a Rajchman measure is not arbitrary and has special porosity. The rate of decay for its Fourier transform is relevant in this respect. See \cite{Ly1} for a background.
 
The Hausdorff dimension of the exceptional set $W(S)$ of  those $x\in\T$ such that $(s_{n}x)$ is not equidistributed and the existence of Rajchman measure supported by $W(S)$ holds our attention in Subsection 2.2.
 \medskip

The above property (\ref{herweyl}) is equivalent to the following one:
\begin{equation} \label{Riemann} \frac{1}{N}\sum_{n=1}^N f(\lambda_{n} x)\to \int_{\T} fdm \quad     a.e.
\end{equation}
for every Riemann-integrable function $f$. But surprisingly, Marstrand (\cite{Ma}), refuting a conjecture of Khinchin, proved that  for $E=\{\lambda_n\}=\N$, there are functions $f\in L^{\infty}(\T)$ for which (\ref{Riemann}) fails. This led us to coin the {\it Khintchin class} of a subset $E=\{\lambda_n\} \subset \N$ as 
\begin{equation}\label{khikhi} \mathcal{K}_E:=\Big\{f\in L^{1}(\T) :  \frac{1}{N}\sum_{n=1}^N f(\lambda_{n} x)\to \int_{\T} fdm \hbox{\quad a.e.} \Big\}
\end{equation} 
In  terms of this notation, the above result of Marstrand can be restated as $\mathcal{K}_{\mathbb{N}}\not\supset L^{\infty}(\T)$. But for the set $S$,   Marstrand (\cite{Ma}) proved that $\mathcal{K}_S\supset L^{\infty}(\T)$ and later Nair (\cite{Na}) proved that 
$ \mathcal{K}_S=L^{1}(\T)$. 
 We will revisit these results in Subsection 2.2.3 and give a simple proof of Nair's result in a weaker form. \\

\noindent $\bullet$ {\bf Harmonic analysis properties}.
 Now consider $S$ as a subset of $\Z$, the dual group of the unit circle $\T$. To better understand  the involved notions, it is useful to consider more generally (see \cite{Ru1},  notably chapter 5, for a detailed exposition)   a compact abelian group $G$ equipped with its normalized Haar measure $m$, its dual $\Gamma=\widehat{G}$, the Fourier transform $\widehat{f}:\Gamma\to \C$ of a function $f\in L^{1}(G)$ or  $\widehat{\mu}:\Gamma\to \C$ of a bounded measure $\mu$ being defined as 
$$\widehat{f}(\gamma):=\int_{G} f(x)\gamma(-x)dm(x),\  \widehat{\mu}(\gamma):=\int_{G} \gamma(-x)d\mu(x).$$
The spectrum of $f$ is (with a similar definition for  $\mu$)  
$$\hbox{\ sp} (f):=\{\gamma\in \Gamma : \widehat{f}(\gamma)\neq 0\}.$$
 For a subset $E$ of $\Gamma$, and a Banach space $X\subset L^{1}(G)$, $X_E$ denotes  the closed subspace  of  $X$ defined by
$$X_{E}=\{f\in X : \hbox{\ sp}(f) \subset E\}.$$
 This subset $E$ is declared ``sparse'',  in  a vague sense depending on $X$, if any $f\in X_{E}$  behaves better than a generic function of $X$. 
 More specifically:
 \begin{itemize}
 \item[(a)] $E$ is called a {\em $\Lambda(p)$-set} if $L^{1}_E\subset L^p$ for some fixed $p>1$.
  \item[(b)] 
$E$ is  called a {\em Sidon set } (or $1$-Sidon set) if  $f\in L^{\infty}_{E}(G)$ implies $\widehat{f}\in \ell^{1}(\Gamma)$.
 \item[(c)] 
 $E$ is  called a {\em $p$-Sidon set} if $f\in L^{\infty}_{E}(G)$ implies $\widehat{f}\in \ell^{p}(\Gamma) 
                 \  (1\leq p<2)$.
  \item[(d)] 
  $E$ is  called a {\em $p$-Rider set} if  $f\in X_{E}(G)$ implies $\widehat{f}\in \ell^{p}(\Gamma)$ 
 where $X$ is the space of randomly bounded functions
 $(1\le p<2)$.
 \end{itemize}
See \cite{Ru1},\ \cite{Ru2} for the first two notions, \cite{Wo},\ (\cite{Bl} p.181) for the third, and \cite{Ro1},\  \cite{LQR},\ \cite{LQR2} for the fourth. 
\medskip

For example, the set $E=\{2^j, j\ge 0\}\subset \Z$ is Sidon  and hence $\Lambda(p)$ for all $p<\infty$. The set $F=\{(2^j, 3^k)\}\subset \Z^2$ is $p$-Sidon exactly for $p\geq 4/3$ (\cite{Bl}). The notion of $p$-Riderness 
coincides with  that of $p$-Sidonicity when $p=1$ \cite{LQ} and is  weaker when $1<p<2$, but reveals easier to handle with (see \cite{Ro1} for a full arithmetic characterization). For example, it provides a simple proof of the non-trivial fact that the union of two Sidon sets is again Sidon \cite{Dr}.\\

Such "sparse" sets $E$ are lacunary in the sense that 
$$I_{E}(N):=|E_N|, \quad {\rm with} \ \  E_N : =E\cap [1, N]$$ is much smaller than $N$. Indeed, we have
\begin{itemize}
\item $I_{E}(N)\ll N^{2/p}$ if $E$ is a $\Lambda(p)$-set ($p>2$).
\item  $I_{E}(N)\ll \log N$ if $E$ is a Sidon  set.
\item $I_{E}(N)\ll (\log N)^{p/(2-p)}$ if $E$ is a $p$-Rider  set ($1\le p<2$).
\end{itemize}
Here is another kind of lacunary sets. We say that $E\subset\Z$ is a {\it Rajchman set} if for any probability measure $\mu$ on $\T$, $\lim_{|n|\to\infty\atop n\in E^c}\widehat\mu(n)=0$ implies that  $\mu$ is a Rajchman measure. 
The relations  {between} these lacunary sets are shown in Figure 1 and the last inclusion is due to Lef\`evre and Rodriguez-Piazza \cite{LR}.
\medskip

For the Furstenberg set $S$, we have $ I_{S}(N) \approx (\log N)^2$, a property  shared by  any  $\frac{4}{3}$-Rider set (NB. $p/(2-p)=2$ for  $p=4/3$). 
This leads naturally to the questions about $\Lambda(p)$-property or $p$-Sidonicity, $p$-Riderness for the set $S$ and its random brother $T$,  which is  defined below. 
We will examine these properties for both $S$ and $T$.
\bigskip

\begin{figure}[H]
\centering
 \begin{tikzpicture}
\begin{scope}

\draw (0,0) circle (3.35);
\draw (0,2.75) node  {Rajchman};
\draw (0,-0.5) circle (2.8);
\draw (0,1.7) node  {$p$-Rider};
\draw (0,-0.95) circle (2.2);
\draw (0,0.65)node{$p$-Sidon};
\draw (0,-1.25) circle (1.5);
\draw (0,-0.2)node{Sidon};
\draw (0,-1.5) circle (1);
\draw (0,-1.35)node{Hadamard};

\end{scope}

\end{tikzpicture}

\bigskip

\caption{Harmonic classification of lacunary sets} 

\end{figure}
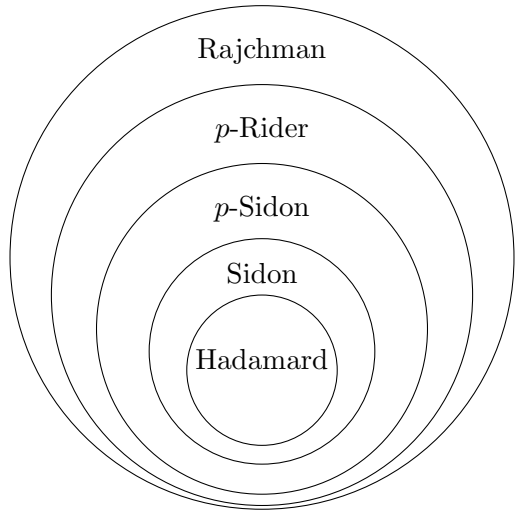

\subsection{The set $T$}
Motivated by the estimate 
$I_{S}(N)\approx (\log N)^2$ seen above, we define  a random version $T$ of $S$ as follows in Section 5: we select the integer $n\geq 1$ with probability $\delta_n=\frac{\log n}{n}$ and we reject it with probability  $1-\delta_n,$ our selections or rejections being independent. This approach is not new. It goes back to Erd\"os-R\'enyi who introduced random sets of integers in \cite{ER}, it  was then systematically developed by Bourgain \cite{Bo},  and  became popular under the name ``Selectors of Bourgain".  See also \cite{LQR}, and the papers  by Kahane and Katznelson \cite{KK1,KK2}. 
We obtain in this way a random set $T$ of integers satisfying
$$\mathbb{E}(I_{T}(N))=\sum_{n=1}^N \frac{\log n}{n}\approx (\log N)^2.$$
It can be proved that $I_{T}(N) \approx (\log N)^2$ almost surely (this  follows from a non-standard law of large numbers (see \cite{BF2005} p.276  for a proof). This random set 
$T$  
 might be easier to deal with, and  might give us some hint on what we could, or could not, hope for the set $S$.
All the properties we have  previously  discussed  will be revisited with $T$ in place of $S$.

\subsection{Detailed content of the paper}
Let us now be more specific on the content of this paper, which is organized as follows.
\medskip

\noi {\bf 1.}  Section 1 is this introduction.
\medskip

\noi {\bf 2.}  Section 2 investigates the first arithmetic or dynamical properties of $S$. 

We naturally start this section with an estimate of the cumulative function  $I_{S}(N)$, 
leading to the asymptotic behaviours of $s_n$, $s_{n+1}/s_n$ and then  of $s_{n+1}-s_n$ and $(s_{n+1}-s_n)/s_n$, improving previous inequalities due to Tijdeman \cite{Ti1,Ti2} by providing effective constants. We obtain the following estimates (Corollary  \ref{TiTi}).
\begin{theo}\label{RMT} We have 
$${1\over(\log s_n)^{\rho}}\ll{s_{n+1}-s_n\over s_n}\ll{1\over(\log s_n)^{1/(\rho+1)}}$$
where the exponent $\rho$ is explicit ($\rho=4.117..$).
\end{theo}
Our estimations relie on a rational approximation of $\alpha=\log2/\log3$ due to G. Rhin \cite{Rh} and Wu and Wang \cite{WW2014}.

Then we switch to the study of the $S$-orbits, that is the sets $(s_nx)$ for $x\in\T$, which is relevant to the dynamics
of the semi-group $\{2^j 3^k\}$. All infinite $S$-orbits $(s_nx)$ for $x\not\in\Q$ are dense in $\T$ but not necessarily equidistributed
(see \cite{Fu1}).  We complete these assertions of Furstenberg by investigating in detail the set $W(S)$ of  those $x\in\T$ such that $(s_{n}x)$ is not equidistributed. As observed above, $W(S)$ is negligible but it is uncountable and even has a positive Hausdorff dimension, for which we give estimates in Section 2 (Theorem \ref{fan}). 
Moreover, we construct a probability measure supported on $W(S)$ with its Fourier transform vanishing at infinity (Theorem \ref{rajmea}).

\begin{theo}\label{baddis} The set of $x$ with a non-equidistributed $S$-orbit has a positive Hausdorff dimension $\ge 0.451621$ and supports a \textnormal{Rajchman} probability measure.
\end{theo}

 Another question concerning the distribution of the orbits $(s_n x)$, related to a conjecture of Khinchin (\cite{Kh}), has been solved in the nineties by Nair \cite{Na}.  Recall that this conjecture, refuted by Marstrand, asked whether $\frac{1}{N}\sum_{n=1}^N f( nx)$ $\to \int fdm \hbox{\ m-a.e.}$  for all $f\in L^1(\T)$. For subsets of $\mathbb{N}$, we can  
speak of the problem of Khintchin. Nair's result solved this problem positively  for the set $S$.
 \begin{theo}[Nair] \label{Nair} For every $f\in L^1(\T)$, 
$\frac{1}{N}\sum_{n=1}^N f(s_nx)$ $\to \int f dm\hbox{\ m-a.e.}$
\end{theo} 
We will give a very simple proof of this fact with a slight restriction on $f$ (membership in  
 $L\log^+ L$) and a release on the measure (Theorem \ref{nair}).

\medskip

\noi {\bf 3.} Section 3 begins with reminders of notions of thin sets from harmonic analysis. 
We recall here different notions of lacunarity: $\Lambda(p)$, $q$-Paley,  $p$-Sidon and later $p$-Rider. The set $S$, thanks to its  weakly lacunary behaviour, is expected to share some of them. Actually, Gundy and Varopoulos \cite{GV} proved in the seventies that $S$ is a $\Lambda(p)$-set  for every $p>2$ (and even $q$-Paley for every $1\le q < 2$) by a judicious interplay of the underlying martingales, and here probability comes into the matter.
\medskip

\noi {\bf 4.} Section 4 presents two probabilistic tools on which we will heavily rely, namely  the Azuma inequalities for scalar martingales with bounded differences, which allow a unified presentation of some results (like Salem-Zygmund's theorem). 
  A second probabilistic tool, the Burkholder inequalities on the square function of a scalar or Hilbert space-valued martingale, is recalled. Then, we prove in  details  the result of  Gundy-Varopoulos  \cite{GV} which deserves to be better known, the initial proof being rather abrupt.
\begin{theo}[Gundy-Varopoulos]\label{rini}The set $S$ is $\Lambda(p)$ for all $2<p<\infty$ and even $q$-Paley for all $1<q<2$. 
\end{theo}
We are not able to prove a $p$-Sidon or $p$-Rider property for $S$, although necessary conditions are shown to be satisfied for an interval of values of the parameter $p$. Whence our last section 5.
\medskip

\noi {\bf 5.} In order to determine which is prevailing between arithmetic or density, we turn to  a random version $T:=T(\omega)=\{t_n\}$ of the set $S$, which a.s. shares with $S$ the same cumulative function; in particular, $t_{n+1}/t_n\to1$ a.s. and almost sure estimates of its growth can be performed (Theorem \ref{diviav} and Theorem \ref{tijdeman}).  Here are some results concerning this set $T$.
\begin{theo}[Bourgain]\label{Bourg} The random set $T=\{t_1<t_2<\cdots\}$ is  a.s. Hartman uniformly distributed i.e. a.s. $\lim_{N\to\infty}{1\over N}\sum_{n<N}e(t_nx)=0$ for every $x\not=0$.
\end{theo}

 This  result, due to Bourgain \cite{Bo}, then revisited in \cite{Ne} and \cite{LQR},  is confirming that the arithmetic property must be relevant in distribution problems since $S$ itself is NOT Hartman uniformly distributed, as proved in Section 2.
 
\begin{theo}\label{bobototo} The random set $T$ admits a.s. a large Khinchin class $X$ containing all $L^p(\T)$ ($p>1$): a.s. $\lim_{N\to\infty}{1\over N}\sum_{n<N}f(t_nx)=0$ a.e. when $f\in X$ satisfies the Khinchin conjecture for the positive integers.
\end{theo}

 This second result, weaker than its deterministic analogue (Nair), seems new and could possibly be improved.  
\begin{theo}\label{llqr} The random set $T$ is a.s. a $p$-Rider set for $p>4/3$.
\end{theo}

 The  above  first two results are valid for more general random sets that we also consider in Section 5. The third result, not shared by all random sets (see \cite{LQR}), acts as an incentive for pursuing with $S$ itself.
On the other hand, $S$ is a  $\Lambda(p)$-set for all $p<\infty$ (see Section 4, Theorem \ref{bernadette}), but we do not know the answer as  concerns $T$. However, it can be proved that $T$ contains a relatively big subset
which is a $\Lambda(p)$-set  for all $p<\infty$ (We will come back to this question in a forth coming paper). 

\section{First results on the Furstenberg set}
The Furstenberg set, denoted by $S$ in this paper,  belongs to a large class of extensively studied multiplicative semi-groups of positive integers. As agreed
in the introduction, we denote by $S(q_1,\dots, q_s)$  the semi-group generated by  some fixed coprime integers $q_1,\dots, q_s$, and by $\{u_1<u_2<\cdots\}$  the sequence 
 obtained from the set $S(q_1,\dots, q_s)$ arranged in increasing order. So, the Furstenberg set $S=S(2,3)$ is the special case of $s=2$ and  $q_1=2, q_2=3$. \\
  Let $p_n$ be the $n$-th prime number.  We know relatively few properties of the set 
 $$S(p_1,\dots, p_s)=:\{u_1<u_2<\cdots\}.$$  
 
 The first observation that $u_{n+1}-u_n\to\infty$ is due to Axel Thue  \cite{Th}, as a consequence of his work on algebraic numbers. Successive improvements have been obtained till the results of Tijdeman \cite{Ti1,Ti2}, which state that
\begin{equation}\label{Ti}{u_n\over(\log u_n)^A}<u_{n+1}-u_n<{u_n\over(\log u_n)^B}
\end{equation}
holds with computable positive constants $A,B$ (for $n$ large enough). Thanks to a geometric consideration, Marstrand \cite[p.545]{Ma}  gave  the following  asymptotic estimate 
  \begin{equation}\label{eq:Marstrand-Est}
 |S(q_1,\dots, q_s)\cap[1,\dots,n]| = K_s (\log n)^s  +r_n \quad
 {\rm as } \ \  n\to\infty,
   \end{equation}
where 
$$
      K_s = \frac{1}{s!} \prod_{j=1}^s\frac{1}{\log q_j},\  r_n=O((\log n)^{s-1}).
 $$
 Marstrand stated the result for $q_j=p_j$, but his proof remains the same for general relatively coprime integers $q_1, \cdots, q_s$.
  
We will propose a more analytic approach to improve the remainder term $r_n$ in (\ref{eq:Marstrand-Est}) in our case $S:=S(2,3)$,  and then deduce better constants  $A,B$ in (\ref{Ti}) for $S(2,3)$.

Dynamics studies the behavior of orbits under a transformation or a semi-group of transformations. 
We study here the action on $\mathbb{T}$ of the semi-group $S$  generated by
$x \mapsto 2x \mod 1$ and $x \mapsto 3x \mod 1$. The orbit of $x$ under $S$ is the set $\{s_n x \mod 1\}$. More generally, 
given a subsequence   of  integers $(k_n)$,  we are interested in  the orbit of $x\in\T$ defined by $(k_nx\mod 1, {n\ge1})$. 
In the case of a Hadamard set, namely $k_n=q^n$ with $q\ge2$, 
the orbits $O_q(x):=(\{q^nx\})$ can be described through the $q$-adic expansion of $x$; it is thus easy to construct uncountably many $x$ such that $O_q(x)$ is not dense. Much better, for example when $q=2$, there exist infinite closed orbits contained in a half-circle, called Sturmian orbits (cf. \cite{BS}).
The other way round, almost all $x$ give rise to a uniformly distributed orbit and the residual set $W(q)$ has Hausdorff dimension 1. 
 Such questions are asked and some results are  obtained in the weakly lacunary case described by the condition $\lim_{n\to \infty}(k_{n+1}-k_n)=+\infty$ (cf. \cite{Ne}).
Furstenberg  \cite{Fu1} was interested in the distribution of the orbits $(s_nx)$, and showed that $(s_nx)$ is dense mod 1 for every irrational number $x$, a first notable difference from the Hadamard case.
We will prove that the set $W(S)$ of $x$ such that $(s_n x)$ is not uniformly distributed is not contained in $\Q$ (in particular, $S=(s_n)$ is not Hartman uniformly distributed) and even has positive Hausdorff dimension.  As we will show, if $S$ were Hartman uniformly distributed, the famous  $\times2\times 3$ Furstenberg  conjecture would be true. 
On the other hand, we will prove that $W(S)$ is a $M_0$-set, meaning that $W(S)$ supports a probability measure whose Fourier coefficients tend to zero at infinity. 

Related to the uniform distribution of the orbits $(s_nx)$ is a question of Khintchin. At the end of this section, we will give a simple proof of the fact, due to  Marstrand \cite{Ma},   that bounded measurable functions are in the Khintchin class  $\mathcal{K}_S$ of the Furstenberg set (see (\ref{khikhi}) for the definition of $\mathcal{K}_S$).

\subsection{ 
Asymptotic properties of the Furstenberg set $S=\{s_n\}$}

\medskip
Throughout this section, we put $\alpha=\log2/\log3$. By a result of Gelfond (\cite{Ge}, Th.~12.2.2  \  p.~226), $\alpha$ is a transcendental number, but not Liouville (or else diophantine).  Recall that a real number $x$ is called $\rho$-diophantine, with $\rho>0$, if 
 for all $n\geq 1$ and a constant $c>0$:
  $$\Vert nx\Vert\geq cn^{-\rho}. $$ 
Any such $\rho$ is called an irrationality exponent of $x$.  We will need the following sharp estimate due to G. Rhin (\cite{Rh} p.~160 with 
$\rho=7.616$), which is improved by Wu and Wang (\cite{WW2014}
with $\rho=4.117$). The first  study on logarithms of integers  is due to  A. Baker \cite{Baker1964}. 

\begin{theo}[\cite{Rh,WW2014}]\label{Rh}
 The transcendental number $\alpha$ is $\rho$-diophantine,
 with $\rho=4.117<\infty.$  
   \end{theo}

Let $S_A:=S\cap [1,\cdots,A]$ ($A$ being an integer). 
According to  Hardy (\cite{Ha} p.~69, notably the equation (5.4.1) there), Ramanujan made the claim 
  $$|S_A|= \frac{\log (2A)\log(3A)}{2\log 2\log 3}+\cdots={1\over 2\log2\log3}\log^2A + 
 {\log6\over 2\log2\log3}\log A +\cdots$$ (where $\cdots$ is a remainder term)  in his first letter to Hardy; the latter obtained $o(\log A/\log\log A)$ for the error term (cf. \cite{Ha}, p.74) and this result stimulated  the study of asymptotic estimates for the cumulative function $|S_A|$.    
We will prove

 \begin{theo}\label{S_N}  Let $\rho$ be an irrationality exponent of $\alpha$, and  $\delta:=\frac{\rho}{\rho +1}$. As $A\to \infty$, we have 
\begin{equation}\label{eq:FQQ}
|S_A|= {1\over 2\log2\log3}\log^2A + 
 {\log6\over 2\log2\log3}\log A +O((\log A)^\delta).
 \end{equation}
 In particular, using Theorem \ref{Rh}, we can take $\rho=4.117$ and $\delta=0.80457299.$
 
 \end{theo}
This result improves  Hardy's error term and the estimate (\ref{eq:Marstrand-Est}) in the case of Furstenberg set $S$.

 Our proof of (\ref{eq:FQQ}) holds for $S(q_1, q_2)$ and the exponent $\delta$ in the error term $O((\log A)^\delta)$ depends on the  irrationality exponent  $\rho$ of $\frac{\log q_1}{\log q_2}$.
 
 \bigskip

   From (\ref{eq:FQQ}), we will deduce the following asymptotic expression for the $n$-th term $s_n$ in $S$. 
 
\begin{theo}\label{BBTA} As $n\to\infty$, 
\begin{equation}\label{nul}  s_n={1\over\sqrt6}\exp(C\sqrt{n}+\delta_n), \hbox{\ where}\  C=\sqrt{2\log2\log3}=1.5230\cdots,
\end{equation}
and  the error term $\delta_n=O(n^{-r})$ with $r={1\over2(\rho+1)}=0.0977\cdots$.
\end{theo}

 In particular, this result recovers a non-lacunarity property of $S$,  namely $s_{n+1}/s_n \to 1$, which was exploited by  Furstenberg \cite{Fu1}. 
  We will finally deduce the inequalities (\ref{Ti}) obtained by Tijdeman, with explicit exponents, in the Furstenberg case.
\begin{theo}\label{TiTi} We have that \textnormal{(}with $2r=1/(\rho+1)$\textnormal{)}
\begin{equation}\label{eq:TiTi}
{1\over(\log s_n)^{\rho}}\ll{s_{n+1}-s_n\over s_n} 
\ll{1\over(\log s_n)^{2r}}
\end{equation}
for $n$ large enough. Recall that $\rho\sim4.117$ and $2r\sim0.1954$.
\end{theo}

\medskip

To prove Theorem \ref{S_N}, we need a Koksma's inequality and an Erd\"{o}s-Tur\'an's inequality.
These two inequalities involve the notion of discrepancy. Let $x_1, \cdots, x_N$ be  a finite set of  numbers in the interval $[0,1]$. Its discrepancy is defined by
$$
     D_N =D_N(x_1, \cdots, x_N) = \sup_{0\le a<b\le 1} \left|\frac{\#\{1\le n\le N: x_n \in [a,b]\}}{N} - (b-a)\right|.
$$
We could say that the supremum is taken over intervals $[a, b]$. If we take the supremum over intervals $[0,b]$, we get a variant $D_N^*$ of $D_N$. 
It is easy to see that $D_N^*\le D_N\le 2D_N^*$ (cf. \cite{KN} p.91).

\begin{lem}
[Koksma inequality, \cite{KN} p.143] Let $f$ be a function defined on the interval $[0,1]$ of bounded variation $V(f)$. Let 
$x_1, \dots, x_N$ be $N$ given points in $[0,1]$ of discrepancy $D_N^*$. Then
\begin{equation}\label{Ko}\Big\vert \frac{1}{N} \sum_{j=1}^N f(x_j)-\int_{0}^1 f(t)dt\Big\vert \leq V(f)D_N^*.
\end{equation}
\end{lem}

Let us state the following special case of Erd\"{o}s-Tur\'an inequality, which says that it is possible to estimate the discrepancy $D_N$ by the Weyl sums. 
\begin{lem}
[Erd\"{o}s-Tur\'an inequality, \cite{KN} p.114] There exists a constant $C>0$ such that for any  $N$ given points
$x_1, \dots, x_N$  in $[0,1]$ and  for any positive integer $m$ we have
\begin{equation}\label{ET}
D_N \le C\left( \frac{1}{m} +\sum_{h=1}^m \frac{1}{h}\left| \frac{1}{N}\sum_{n=1}^N e^{2\pi i h x_n}\right|\right).
\end{equation}
\end{lem}

{\em Proof of Theorem \ref{S_N}.} 
 Denote
$N_A:=
\#\{(j,k)\ge0 \ ;\  2^j3^k\le A\}$. Clearly 
$$N_A=\#\{(j,k)\ge0 \ ; \  j\log 2+k\log 3\le \log A\}
=\#\{(j,k)\ge0\ ; \  j\alpha+k\le B\alpha\},$$
with  $B=\log A/\log 2$. Hence we get
$$N_A=\sum_{k+j\alpha\le B\alpha}1=\sum_{0\leq j\le B}\big([(B-j)\alpha]+1\big).$$
Using the relation $[x] = x -\{x\}$, we get
\begin{equation}\label{bel}
N_A=\sum_{0\leq j\le B}\big(1+(B-j)\alpha\big)-\sum_{0\le j\le B}\{(B-j)\alpha\}.
\end{equation}
\noi 
Let us first estimate the second sum on the right hand side of (\ref{bel}), the first sum being easy to compute. Write 
$$\sum_{0\le j\le B}\{(B-j)\alpha\} =
 \sum_{j=1}^N f(x_j)+O(1)$$ with $f(x)=\{x\}$,
$x_j=(B-j)\alpha$ and $N=[B]$. We claim that
 
as  $B\to \infty$, we have
  \begin{equation} \label{top} \sum_{0\le j\leq B} \{(B-j)\alpha\}=\frac{B}{2}+O(B^\delta), \ \delta={\rho\over\rho+1}\cdot
  \end{equation} 
Indeed, as $\int_0^1f(x)dx =\frac{1}{2}$ and $V(f)=1$,
by the
Koksma inequality (\ref{Ko}), we get 
\begin{equation}\label{Ko2}\Big\vert \frac{1}{N} \sum_{j=1}^N f(x_j)-\frac{B}{2}\Big\vert \leq  D_N;
\end{equation}
 
Now, by the Erd\"{o}s-Tur\'an inequality (\ref{ET}), the discrepancy  $D_N$ of the finite sequence $x_j=(B-j)\alpha \mod 1$ ($1\leq j\leq N)$ is, for every $m\in \N$, bounded  above by 
   $$D_N 
   \ll  \frac{1}{m}+\sum_{h=1}^m \frac{1}{N\,h\Vert h\alpha\Vert}$$
   (cf. \cite{KN} p.123);
but  Theorem \ref{Rh} implies $\sum_{h=1}^m 1/{h\Vert h\alpha\Vert}\ll\sum_{h=1}^mh^{\rho-1}\ll m^\rho$ which provides
   $$D_N\ll  \frac{1}{m}+\frac{m^{\rho}}{N}\cdot $$ 
 We optimize this quantity by taking $m=N^{1/(\rho+1)}$ and obtain 
   $D_N\ll N^{-1/(\rho+1)}$;  we conclude (\ref{top}) with the aid of  (\ref{Ko2}).

 \smallskip   
Back to (\ref{bel}). By using (\ref{top}) we get
$$N_A=B+ \alpha\Big(B(N+1)-\frac{N(N+1)}{2}\Big)-{B\over2}+O(B^\delta)$$
whence
\begin{equation}\label{soen}N_A= \frac{\alpha+1}{2}B+\alpha\big(BN-\frac{N^2}{2}\big)+O(B^\delta);
\end{equation}
now, writing
$BN-N^2/2=-{1\over2}\big[(N-B)^2-B^2\big]=B^2/2+O(1),$
we deduce from (\ref{soen}) the asymptotic behaviour:
\begin{equation} \label{noen}N_A=aB^2+bB+O(B^\delta),
\end{equation}
with $B=\log A/\log 2$, $\delta={\rho\over\rho+1}$, $a=\alpha/2$ and $b=(1+\alpha)/2$. This proves Theorem  \ref{S_N}.
\hfill$\square$
\medskip

{\em Proof of Theorem \ref{BBTA}.}
We deduce  Theorem \ref{BBTA} from Theorem \ref{S_N}. Take $A=s_n$; thus $B=B_n=\log s_n/\log2$ and $N_A=n$. The  equality (\ref{noen}) gives first $B_n=O(\sqrt{n})$ and then
$$n=aB_{n}^2+bB_n+\alpha_n \hbox{\ with}\ \alpha_n= O(n^{\delta/2}).$$
Consider $B_n$ as unknown of a  second degree algebraic equation.
Solving  $B_n$ gives
$B_n=\frac{1}{2a} \big(-b+\sqrt{\Delta_n}\big)$
with 
$$\sqrt{\Delta_n}=\sqrt{b^2 +4a(n-\alpha_n)}=2\sqrt{a n}\sqrt{1+\frac{b^2}{4an}-\frac{\alpha_n}{n}}$$$$=2\sqrt{a n}\Big(1-\frac{\alpha_n}{2n}+O(n^{-1})\Big)=2\sqrt{a n}+O(n^{-\frac{1-\delta}{2}}).$$
 Here we have used the fact that $2\delta >1$. We then get
 \begin{equation}\label{neige}B_n=\sqrt{\frac{n}{a}}-\frac{b}{2a}+\delta_n \hbox{\ with}\ \delta_n=O(n^{-\frac{1-\delta}{2}})\end{equation} 
This gives Theorem \ref{BBTA}, since $\frac{1-\delta}{2}=r$ and $B_n=\log s_n/\log2$.
$\square$
\medskip

{\em Proof of Theorem \ref{TiTi}.} 
$\triangleright$ We easily deduce from Theorem \ref{BBTA} the right-hand inequality in (\ref{eq:TiTi}). Since
$${s_{n+1}-s_n\over s_n}\ll C(\sqrt{n+1}-\sqrt n)+|\delta_{n+1}|+|\delta_n|$$
with $\delta_n=O({1\over n^r})$ and $r=\frac{1}{2(\rho+1)}\in(1/18, 1/17)$. With no more information on the sign of $\delta_n$,  we get at best  the estimate  $\frac{O(1)}{ n^r}$.

\smallskip
$\triangleright$ For the reverse inequality, we 
 write $s_n=2^{c_n}3^{d_n}$ and $s_{n+1}=2^{c_{n+1}}3^{d_{n+1}}$ and observe that $c_{n+1}-c_n\not=0$ (which is indeed clear: if $c_{n+1}=c_n$, then $d_{n+1}\geq d_n +1$ and $s_{n+1}\geq 3 s_n$, while clearly $s_{n+1}\leq 2s_n$ for all $n$, since $2s_n\in S$), obviously with the same sign as $d_n-d_{n+1}\not=0$.
 Now, we bound
$${s_{n+1}-s_{n}\over s_n}>\log{s_{n+1}\over s_n}=\big|(c_{n+1}-c_n)\alpha-(d_n-d_{n+1})\big|\cdot \log 3$$
$$=\Big||c_{n+1}-c_n|\alpha-|d_n-d_{n+1}|\Big|\cdot \log 3\ge\|(c_{n+1}-c_n)\alpha\| \cdot \log 3 \ge c|c_{n+1}-c_n|^{-\rho}.$$
(Once again, we make use of G.~Rhin's or Wang-Wu's lower bounds).
Since
$$|c_{n+1}-c_n|\le\max(c_n,c_{n+1})\le\max(\log_2s_n,\log_2s_{n+1})$$
we obtain
$$\|(c_{n+1}-c_n)\alpha\|\ge C(\log s_{n+1})^{-\rho}\sim C(\log s_{n})^{-\rho}. $$
This provides the left-hand inequality with exponent $\rho$; in particular, we deduce $s_{n+1}-s_n\to\infty$.
$\square$
\medskip

\noi{\bf Remark.} We can add a small precision to (\ref{TiTi}): for infinitely many pairs $(s_n,s_{n+1})$, we have
$${s_{n+1}-s_n\over s_n}\le {2\log2\log3\over\log s_n}\cdot$$
The proof goes as follows: due to the best approximation property of the convergents to $\alpha$, one easily sees that $2^q, 3^p$ or $3^p,2^q$ are consecutive terms in $S$ 
for any convergent $p/q$ to $\alpha$ (according to the parity of the indices) (\cite{BDH}). Now, assume for example $s_n=2^q<3^p=s_{n+1}$ to be consecutive in $S$.
Since $|\alpha-p/q|<1/q^2$, we get
$$\log{s_{n+1}\over s_n}=\log{3^p\over 2^q}=\log3\times |q\alpha-p|<{\log3\over q}={\log2\log3\over\log s_n}\cdot$$
Hence the result follows since ${s_{n+1}-s_n\over s_n}\le2 \log{s_{n+1}\over s_n}$.

\subsection{Hartman uniform distribution and Furstenberg conjecture}
It is amazing to notice that the Furstenberg conjecture would be implied by the assertion that the Furstenberg set is Hartman uniformly distributed. But the Furstenberg set is not Hartman uniformly distributed, as we 
will see.

We start with a few reminders on a notion of distribution for {\it sequences of integers}, studied by Hartman.
Let $E:=(k_n)\subset\Z$ and $E_N=E \cap[1,\dots,N]$. Recall that $\T$ denotes the set of reals modulo one with $\Z$ as its dual group.
\begin{defi} We say that the sequence $E$ is Hartman uniformly distributed ($H$-ud in short) if one of the following equivalent conditions is satisfied:\\
\indent \textnormal{ (i)} For every $x\in\T\setminus \{0\}$, ${1\over N}\sum_{n=1}^Ne(k_nx)\to 0$ when $N\to\infty$.\\
\indent \textnormal{ (ii)} The probability measure ${1\over N}\sum_{n\le N}\delta_{k_n}$ is  w$^*$-convergent  to the Haar measure 
 of the Bohr compactification $\overline{\Z}$ of the group $\Z$.
 \end{defi}
In particular, such a sequence is dense in the Bohr compactification of $\Z$, the dual group of $\T$ equipped with the discrete topology. We will give more details about the  Bohr topology in the last section of the paper.
\medskip

\noi{\bf Example 1.}  Of course $E=\N$ is $H$-ud; but the set of squares is not.  Take $x=1/4$ and note that ${1\over N}\sum_{n=1}^Ne(n^{2}x)\to \frac{1+i}{2}$.

\medskip
\noi{\bf Example 2.}  Let us consider the Rudin-Shapiro sequence $(r_n)$, which is defined by $
r_n=(-1)^{u_n}$ where $u_n = \sum_{k\ge 0}\epsilon_k(n)\epsilon_{k+1}(n)$ for $n = \sum_{k\ge 0} \epsilon_k(n)2^k$ (the dyadic expansion of $n$). 
It is known that $\Vert \sum_{n\leq N} r_n e_n\Vert_\infty \ll \sqrt{N}$  (\cite{Ru2}). We choose as $E =(n_k)$ the set of occurrences of 1 in the sequence, i.e. $r_{n_k}=1$, $r_n=-1$ else. This set $E$ is H-ud. Indeed, observe that
$$\sum_{n\le N} r_ne(nx))
=\sum_{n\in E_N} e(nx)\  - \sum_{n\le N\atop n\notin E_N} e(nx)=: S_N'(x)-S_N''(x).$$
In addition, we have $\sum_{n\le N} e(nx) = S_N'(x)+ S_N''(x)$.
It follows that
$$\sum_{n\in E_N} e(nx)={1\over2}\Big(\sum_{n\le N} e(nx)+\sum_{n\le N} r_ne(nx)\Big)
=O({1})+O(\sqrt{N} )
$$
 for $x\neq 0$.
Then we conclude  by using the relation (the above first equality when $x=0$): $2|E_N|=N+\sum_{n\leq N}r_n\approx N$

\medskip
\noi{\bf Example 3.}  The Thue-Morse sequence $(t_n)$ is defined by $
t_n=(-1)^{s(n)}$ where $s(n) = \sum_{k\ge 0}\epsilon_k(n)$ for $n = \sum_{k\ge 0} \epsilon_k(n)2^k$ (the dyadic expansion of $n$). 
It is known that $\Vert \sum_{n\leq N} t_n e_n\Vert_\infty \ll N^{\log 3/\log 4}$  (\cite{Gelfond}). By the same argument as in {\bf Example 2},
we can prove that  the set   $E =(n_k)$ of occurrences of 1 in the sequence $(t_n)$ is H-ud. The result in \cite{FSS} on the supremum norm of trigonometric polynomials with generalized Thue-Morse coefficients can be similarly used to get more H-ud sequences from  generalized Thue-Morse sequences.

\bigskip
\noi{\bf Furstenberg conjecture.}
Making use of the fact $s_{n+1}/s_n\to1$, Furstenberg proved that {\it every closed infinite set of $\T$ which is $S$-invariant must be $\T$ itself} (see \cite{Bos}) (``$F$ is $S$-invariant" meaning $2F\subset F$ and $3F\subset F$). See \cite{BLMV} for effective versions of this result. Of course there are  lots of closed  infinite proper subsets only $2$-invariant (or $3$-invariant) and this deep difference led Furstenberg to his famous metric conjecture: 

 \begin{quotation} \it A continuous probability measure on $\T$ which is $S$-invariant
  must be the Lebesgue measure $m$.
  \end{quotation}
\noi (That $\mu$ is $S$-invariant or $\times2\times3$-invariant means $\widehat\mu(sn)=\widehat\mu(n)$ for every $s\in S$, $n\in\Z$). 
It is interesting to point out a connection between Hartman uniform distribution and this latter conjecture. 

{\it If $S=(s_n)$ were Hartman uniformly distributed, in other terms, if we had
\begin{equation}\label{star}P_{N}(x):=\frac{1}{N}\sum_{n=1}^N e(s_{n}x)\to 0\ \  \hbox{\ for}\ x\in \mathbb{T}\setminus\{0\},
\end{equation}
then the Furstenberg conjecture would be true.
}

Indeed, suppose that $S$ is $H$-ud and that there exists a $\times2\times 3$-invariant continuous probability measure  $\mu$ different from the Lebesgue measure. 
Then there exists a positive integer $a\neq 0$ such that $\widehat{\mu}(a)\neq 0$. Consider the probability measure $\nu=\mu\circ \sigma_{a}^{-1}$, the image of $\mu$ under the transformation $\sigma_{a}:\T\to \T$ defined by  $\sigma_a(x)=ax$ mod 1, which is still continuous. Since  $\widehat{\nu}(n)=\widehat{\mu}(an) \hbox{\ for all}\ n \in \Z$, the measure $\nu$ is 
again $\times2\times 3$-invariant, with $\widehat{\nu}(1)\neq 0$.   Now, we observe that, by the $\times2\times 3$-invariance of $\nu$,
$$\int_{\T}P_{N}(-x)d\nu(x)=\frac{1}{N} \sum_{n=1}^N \widehat{\nu}(s_{n})=\widehat{\nu}(1)\neq 0.$$
While, by (\ref{star}) and Lebesgue's dominated convergence theorem, the LHS tends to 
$\nu(\{0\})$, 
 a contradiction since $\nu$ is continuous.
\medskip

We shall see that the Furstenberg set $S$ is far from $H$-ud which is  a heuristic indication that Furstenberg's conjecture could be difficult to prove.  We enlarge to the Weyl uniform distribution $\bmod \ 1$ of the orbits $(k_n x)$ for some sequence $(k_n)$ of integers. Thanks to the classical Weyl's theorem, this means
$$
\forall h\in \mathbb{Z}\setminus\{0\}, \quad \ {1\over N}\sum_{n<N}e(hk_nx)\to 0\quad  a.e.
$$
Observe the relation between  the Weyl uniform distribution and  the Hartman uniform distribution: {\it  if $(k_n)$ is H-ud, then the sequence $(k_nx)$ is uniformly distributed mod 1 for every irrational $x$.}
The Furstenberg sequence $\{s_n\}$ is not H-ud, because there are uncountably many $x$ such that $\{s_n x\}$ is not  Weyl uniformly distributed, as we  will see below.

\subsection{The set $W(S)$ of non $S$-normal numbers}
We focus on the negligible set of non $S$-normal numbers:
$$
W(S):=\left\{x\in \mathbb{T}\ :\ \{s_nx\} \ {\rm is\ not\ uniformly\ distributed}\right\}
$$ 
At the end of his article (\cite{Fu1}) Furstenberg pointed out that $\{s_nx\}$ is not uniformly distributed mod 1 for some suitable Liouville number. 
As an immediate consequence, 
{\it the set $S$ is not Hartman uniformly distributed}. 
We can say more than the above remark of Furstenberg.
First, by constructing a lot of such Liouville numbers, we observe that the set  $W(S)$
is uncountable; this was already noticed in \cite{Bos} and recently, C.~Badea and S.~Grivaux \cite{BG} recover this fact in a non-explicit way (see also \cite{BGM}). We can go further by examining the Hausdorff dimension of $W(S)$.
In this respect, Erd\"os and Taylor (\cite{ET, Fan}) proved that the negligible set of non $\Lambda$-normal numbers has Hausdorff dimension $1$  when $\Lambda\subset\Z$ is a Hadamard set; this result holds true in many other cases, for example the set $W(\Lambda)$ for $\Lambda=\{3^j+3^k\}_{j,k\ge0}$ has the same dimension 1 although $\Lambda$ is not Hadamard (\cite{Ne}). 
Boshernitzan \cite{Bos} announced the result   $\dim_HW(S)=1$, due to D. Berend. But  D. Berend (personal communication) claims  that he  never thought  he could prove it.

 \begin{theo}\label{fan} The Hausdorff dimension of $W(S)$ is $\ge 0.451621$.
 \end{theo}
 \begin{proof} 
 We are going to 
construct inside $W(S)$ a kind of Cantor set, called a homogeneous Moran set \cite{FWW,FRW}. Here is a general construction (cf. \cite{FRW, FWW}). 
Let $( n_k)_{k\geq 1}$ be a
sequence of positive integers  with $n_k\ge 2$ and $( c_k)_{k\geq 1}$ be a sequence of
positive numbers satisfying $0<c_k<1$
and $n_kc_k\leq 1$ ($k\geq 1$). 
For $k\ge 1$, let $D_k =\{(i_1, \dots, i_k): 1\le i_j \le n_j \ {\rm for} \ 1\le j \le k\}$ and $D=\bigcup_{k\ge 0} D_k$ with $D_0=\emptyset$. 
For $\tau=(\tau_1, \dots, \tau_n) \in D_n$ and  $\sigma=(\sigma_1, \dots, \sigma_m) \in D_m$, we define
$\tau*\sigma =  (\tau_1, \dots, \tau_n, \sigma_{1}, \dots, \sigma_{m}) \in D_{n+m}$. Suppose we are given an interval $J$ of length $1$. We can construct 
a family of subintervals $\mathcal{F} :=\{J_\sigma: \sigma\in D\}$ of $J$ as follows. The construction will be made for   $\{J_\sigma: \sigma\in D_k\}$ inductively on $k$. 
First for $k=0$, we choose $J_\emptyset =J$. Now suppose   that $\{J_\sigma: \sigma\in D_{k-1}\}$ are constructed. For each $\sigma \in D_{k-1}$, we choose
$n_k$ subintervals $J_{\sigma *1}, \dots, J_{\sigma* n_k}$ in $J_\sigma$ such that they have disjoint interiors and 
$$
     \forall 1\le j \le n_k, \ \ \ \frac{J_{\sigma*j}}{|J_\sigma|} = c_k.
$$
The set $E:=\bigcap\limits_{k\geq 1} \bigcup\limits_{{\ \sigma }\in D_k}J_ { \sigma }$
is called a \textsl{homogeneous Moran set} determined by $\mathcal{F}$.
Notice that for given sequence $(n_k)$ and $(c_k)$, there are different constructions of $\mathcal{F}$ and then different Moran sets, because the positions 
$J_{\sigma *1}, \dots, J_{\sigma* n_k}$ in $J_\sigma$ are arbitrary to some extent.
There is a common lower bound for the Hausdorff dimensions of these Moran sets, given by the following proposition.

\begin{lem} [\cite{FRW,FWW}] \label{prop:Moran}
For the homogeneous Moran set defined above,
we have
\[
\dim_H E\geq \liminf_{n \rightarrow \infty} \frac{ \log n_1n_2\cdots n_k}{
-\log c_1c_2\cdots c_{k+1}n_{k+1}}\cdot
\]
 \end{lem}
 
 Now let us construct a Moran set in our set $W(S)$.
Let us fix a  sequence of integers $(\ell_k)$ which is quickly increasing in the sense that 
$$
   \ell_0=0, \qquad \frac{\ell_1+\cdots +\ell_{k-1}}{\ell_k} \to 0.
$$
For a small $\delta>0$, define $m_k = \delta \ell_k$ so that $\ell_{k-1}<m_k<\ell_k$. 
We consider
$$
    E :=\left\{\sum_{j=1}^\infty \frac{a_j}{6^j}:  a_j =0 \  {\rm when} \ j \in (m_{k}, \ell_{k}) \ {\rm for \ some}\ k,\ 0\le a_j\le 5 \ {\rm otherwise}\right\}\!\cdot
$$
The set $E$ is a Moran set with
 $$
   n_k = 6^{m_k -\ell_{k-1}}, \qquad c_k = \frac{1}{6^{m_k -m_{k-1}}}\cdot
$$
Indeed,  in the definition of $x= \sum_{j=1}^\infty {a_j}6^{-j}\in E$, the digits $a_j$ are free (meaning they  can take any value between $0$ and $5$) for $j\in [\ell_{k-1}, m_k]$.
So there are a number $n_k$ of possible strings of digits $(a_j)$ with $ j\in [\ell_{k-1}, m_k]$ corresponding to subintervals of an $(k-1)$-level interval. 
We have then
\begin{eqnarray*}
\frac{\log n_1n_2\cdots n_k }{-\log c_1c_2\cdot c_{k+1} n_{k+1}}
  & = &\frac{\sum_{j=1}^k (m_j -\ell_{j-1})}{ \sum_{j=1}^{k+1} (m_{j} -m_{j-1}) +\ell_{k}-m_{k+1}}\\
  & = &\frac{\sum_{j=1}^k (m_j - \ell_{j-1})}{ \ell_{k}}\cdot
\end{eqnarray*}
By Lemma \ref{prop:Moran}, 
$$
   \dim_{H} E \ge \liminf \frac{\sum_{j=1}^k (m_j - \ell_{j-1})}{ \ell_{k}}=
   \liminf  \frac{\sum_{j=1}^k (\delta \ell_j - \ell_{j-1})}{ \ell_{k}} = \delta.
$$
Next we show that $E\subset W(S)$ if $\delta>0$ is small enough. Let us choose $p_k = (1-\epsilon) \ell_k$
with $\epsilon>0$ small enough to have $m_k <p_k <\ell_k$. Observe that
$$
\ell_k -p_k = \epsilon \ell_k  \to \infty, \qquad \frac{m_k}{p_k} \le \frac{\delta}{1-\epsilon}.
$$
Let 
$x= \sum_{j=1}^\infty {a_j}6^{-j}\in E$.  For every $k$, we have
\begin{equation} \label{alors2}
    x= \sum_{j=1}^{m_k} \frac{a_j}{6^j} +  \sum_{j=\ell_k}^\infty \frac{a_j}{6^j} =:a_k^* 6^{-m_k}+O(6^{-\ell_k})
\end{equation}
where $a_k^*$ is an integer. Then we put $N_k=6^{m_k},\ N=6^{p_k}$ and consider
\begin{equation}\label{split}
T_N=\sum_{2^{i} 3^{j}<N} e(2^{i} 3^{j}x)=:u_N+v_N+W_N-z_N,
\end{equation} 
with
$$u_N=\sum_{i< m_k, 2^{i}3^{j}<N} e(2^{i} 3^{j}x),\  \ v_N=\sum_{j< m_k, 2^{i}3^{j}<N} e(2^{i} 3^{j}x), $$
$$ W_N=\sum_{i,j\geq m_k, 2^{i}3^{j}<N} e(2^{i} 3^{j}x), \quad
z_N=\sum_{i< m_k, j< m_k} e(2^{i} 3^{j}x).$$
It is clear that (notice that $\log N=p_k \log 6$)
\begin{equation}\label{clair2} |u_N|\le {\frac{\log 6}{\log 3}} m_kp_k, \quad\ \  |v_N|\le {\frac{\log 6}{\log 2}}m_k p_k, \quad \ \ z_N\le m_k^2.
\end{equation}
Now, if $i\geq m_k$, $j\geq m_k$ and $2^{i}3^{j}<N$, we get from (\ref{alors2}):
$$ e(2^{i}3^{j}x)=1+O\big(\frac{N}{6^{\ell_{k}}}\big)=1+O(6^{p_{k} -\ell_{k}})=1+o(1), $$
because  $\ell_{k}-p_k\to \infty$. It follows that
\begin{equation}\label{alor}W_N=C\log^{2}N+o(\log^{2}N),\hbox{\ with}\ C=\frac{1}{2\log 2\log 3}.
\end{equation}
 On the other hand, 
 $$
\frac{m_kp_k}{\log^2 N} \le \frac{\delta}{(1-\epsilon) \log^2 6}, \quad 
\frac{m_k^2}{\log^2N} \le \frac{\delta^2}{(1-\epsilon)^2 \log^2 6 }.$$ 
Combining this,  (\ref{clair2}) and (\ref{alor})  results in
$$
\liminf_{k\to \infty} \frac{|T_N|}{\log^{2}N}\geq C- \left({\frac{\log 6}{\log 3} + \frac{\log 6}{\log 2} } \right) \frac{\delta}{(1-\epsilon) \log^2 6} -  \frac{\delta^2}{(1-\epsilon)^2 \log^2 6}>0
$$
if $\delta>0$ is small enough. This proves $x \in W(S)$.  

An explicit lower bound of the dimension can be obtained as the positive  solution of
$$
C-  \frac{\delta}{\log2\log3}-  \frac{\delta^2}{ \log^2 6}=0,
$$
i.e. $ \delta^2 + 2a \delta -a=0$ with $a=\frac{\log^2 6}{2\log 2 \log 3}$,   whence $ \delta =0.451621...$ 
\end{proof}

\medskip
\subsection{$W(S)$ is a $M_0$-set}
 A $M_0$-set is a Borel subset of $\T$ that supports a non-trivial Rajchman measure. Recall that a Rajchman measure (or $M_0$-measure) on $\T$ is a measure $\mu$ such that 
$$\lim_{|n|\to\infty}\widehat\mu(n)=0.$$
The question ``what does the support of a Rajchman measure look like?" was asked by Kahane and Salem \cite[p.59]{KS1964}.
Russel Lyons  \cite{Ly} observed  that a measure annihilating all non-normal sets $W(\Lambda)$ must be a Rajchman measure, but the converse is not true, since a Rajchman measure can be constructed on the set of non-normal numbers to base 2. The construction relies on the 2-adic expansion and the divisibility property of $\{2^n\}$.
By using the 6-adic expansion as above, we can exhibit a Bernoulli-like probability measure which is Rajchman and concentrated on the set $W(S)$.

\begin{theo}\label{rajmea} The set $W(S)$ supports a Rajchman measure, briefly, $W(S)$ is a $M_0$-set.
\end{theo}

\begin{proof} The candidate will be the distribution $\mu:=\P_X$ of a series of discrete random variables
$X:=\sum_{j=1} ^\infty\eps_jX_j$
where $(\eps_j)$ is a sequence of Bernoulli variables such that
$$\P(\eps_j=0)=1/j\log j,\qquad  \P(\eps_j=1)=1-1/j\log j,$$
and the variables $X_j$ are to be defined, all the variables involved being independent.
We fix a lacunary sequence of integers $(m_k)$ so that $m_k/m_{k+1}\to0$ and we put $n_k:=m_{k+1}-m_k$;
the variable $X_ k$ is equidistributed on the finite set
$\{\frac{a}{ 6^{m_{k+1}}}: 0\le a \le 6^{n_k}-1\}$, namely
$$\P(X_k={a\over 6^{m_{k+1}}})={1\over 6^{n_k}}.$$
For every $\omega$, $X(\omega)$ is a 6-adic expansion, with a large block of arbitrary digits.
We have to check that $\mu$ is concentrated on $W(S)$ and that $\mu\in M_0$.

\bigskip
\noi {\bf1.}  We prove that $\limsup{1\over |S_N|}\Big|\sum_{2^i3^j\le N}e(2^i3^jX)\Big|>0$ a.s. which implies that $X$ a.s. takes its values in $W(S)$ (recall that $S_N=S\cap[1,\dots,N]$).\\
 
 As in  (\ref{split}), we write
$\sum_{2^i3^j\le N}e(2^i3^jX)=u_N+v_N+W_N-z_N$
 with now
\begin{equation}\label{clair3} |u_N|\le {\frac{\log N}{\log 3}} m_k, \quad\ \  |v_N|\le {\frac{\log N}{\log 2}}m_k, \quad \ \ |z_N|\le m_k^2.
\end{equation}
 We then diverge from (\ref{split}) and adopt the following notations:
 $$ M_k=6^{m_k},\  N_k=6^{m_{k+1}} 2^{-k-1},\  Q_k=M_{k+1}/M_k= 6^{n_k}.$$
 Since $\sum P(\epsilon_k=0)=\infty$, almost surely $\epsilon_k=0$ for infinitely many $k's$ (say for $\omega\in \Omega_0,\  \mathbb{P}(\Omega_0)=1$). Let now $\omega\in \Omega_0$ and  $k$ large with $\eps_k=\eps_{k}(\omega)=0$.
We take $N=N_k$ in (\ref{clair3}). Once more it holds
$$ e(2^{i}3^{j}X)=e(2^{i} 3^j\sum_{n\ge k+1}\eps_nX_n)\ \ {\rm if}\  i\geq m_k,\ j\geq m_k;$$
and, since $2^{i}3^{j}<N$, we get  
that,  for $n\geq k+1$: 
$$2^{i} 3^j X_{n}\in\Big\{0,{2^{i}3^j\over6^{m_{n+1}}},\ldots,2^{i}3^j\big({1\over6^{m_{n}}}-{1\over6^{m_{n+1}}}\big)\Big\}
$$
with
$$2^{i}3^j\sum_{n\geq k+1} \frac{1}{6^{m_n}} <{2N\over6^{m_{k+1}}}={1\over2^{k}},$$
so that
\begin{equation}\label{mtge}2^{i} 3^j \big(\sum_{n\ge k+1}\eps_nX_n\big)\in[0,1/2^k)\ \ \hbox{for such a}\ k.
\end{equation}
We now focus on the main sum $W_N$.  We begin by observing that\\
 $N=6^{m_{k+1}}2^{-k-1}\in]M_k,M_{k+1}[$,  since $m_{k+1}-m_k\gg k$.  We claim that
$$\limsup_N{1\over |S_N|}\Big|\sum_{i,j\geq m_k, 2^{i}3^{j}<N} e(2^{i} 3^{j}X(\omega))\Big|>0.
$$
Indeed, by (\ref{mtge}),
$$|W_{N_k}|=\Big|\sum_{i,j\geq m_k, 2^{i}3^{j}<N_k} e(2^{i} 3^{j}\sum_{n\ge k+1}\eps_nX_n)\Big|=\sum_{2^{i}3^j<N_k\atop i,j>m_k}1 + O(2^{-k}|S_{N_k}|).$$
Combined with (\ref{clair3}), this gives that, almost surely, for infinitely many $k$'s,
$$\big|\sum_{2^i3^j\le N_k}e(2^i3^jX)\big|=|S_{N_k}|+O(2^{-k}|S_{N_k}|)+O(m_km_{k+1})$$
whence
 $${1\over|S_{N_k}|}\big|\sum_{2^i3^j\le N_k}e(2^i3^jX)\big|=1+O(2^{-k})+O(m_k/m_{k+1})
$$
since $|S_{N_k}|\approx\log^2N_k\approx m_{k+1}^2$. 
The claim follows from $\lim_k m_k/m_{k+1}=0$.

\bigskip
\noi {\bf2.}  We now prove that $|\widehat\mu(h)|\le \phi(h)$ where $\lim\phi(h)=0$ as $h\to \infty$.

For $h\not=0$, 
$$\widehat\mu(h)=\prod_{j=1}^\infty {\E}(\exp 2i\pi h\epsilon_jX_j)=\prod_{j=1}^\infty [(1-1/j\log j)c_j+1/j\log j]$$
where $c_j=\E(\exp 2\pi i h X_j)$. Clearly,
$$c_j={1\over Q_j}\sum_{a=0}^{Q_j-1} e^{2\pi i h{a\over M_{j+1}}}={1\over Q_j}{e^{2\pi i {hQ_j\over M_{j+1}}}-1\over e^{2\pi i {h\over M_{j+1}}}-1}$$
and
$$|c_j| ={M_j|\sin(\pi h/M_j)| \over
M_{j+1}|\sin(\pi h/M_{j+1})|}$$
\noi Assuming $h>0$, we fix the index $k$ such that $M_k/2 < h \le M_{k+1}/2$; then,
$$ 2h> M_k,\  \pi h/M_{k+1} \le \pi/2 $$ 
whence
$$M_{k+1}\sin(\pi h/M_{k+1})\ge 2h$$
by concavity of the sine function on $[0,\pi/2]$.
It ensues that 
$$|c_{k-1}c_k|={M_{k-1}|\sin(\pi h/M_{k-1})|\over M_{k+1}|\sin(\pi h/(M_{k+1})|}\le{M_{k-1}\over2h} \le
{M_{k-1}\over M_k}=:\delta_k;$$ 
this implies that one term in the product $\prod_j [(1-1/j\log j)c_j
+1/j\log j]$ must be small when $h$ becomes large :
indeed, one among both coefficients $|c_k|$ and $|c_{k-1}|$,  say $|c_k|$, must be $\le \sqrt{\delta_k}$ so that, 
$$|\widehat\mu(h)|\le|(1-1/k\log k)c_{k}+1/k\log k|$$
$$\le\sqrt{\delta_k}+{1\over k\log k}=O\Big({1\over k\log k}\Big)=O\Big({1\over\log\log h}\Big)$$
by choosing $m_k=k!$ and using the Stirling formula. The exponent $1$ of the iterated logarithm  is optimal by the Davenport-Erd\"os-Lev\^eque criterion \cite[Chap.1]{Bug 2} 
\end{proof}

It can be more precise about the order of $\frac{1}{\log\log h}$ appearing at the end of the above proof. 

 \begin{theo}\label{upr} The following holds:\\
  1.~  $\widehat{\mu}(h)$ can decay as $1/\log\log h$ as $h\to \infty$.\\
  2.~ $\widehat{\mu}(h)$ cannot decay as $1/(\log\log h)^\alpha$ as $h\to \infty$ as soon as $\alpha>1$.
   \end{theo} 
  \begin{proof} 1. has just been proved: in our construction,  changing our choice $m_k=k!$ to $m_k=\exp(k\log\log k)$ (integral part of), we get in principle 
  $$\widehat{\mu}(h)=O\Big(\frac{1}{\log \log h}\times \frac{\log_{4}h}{\log_{3}h}\Big)$$
  but 2. will show that we  can hardly do better. This is a manifestation of the uncertainty principle: $\widehat{\mu}$ cannot be too small as soon as the support of $\mu$ is porous.
  Note in passing  that, if an "exotic" probability measure $\mu$ exists,  namely $S$-invariant,  $S$-ergodic and $\neq m$,  it must be supported by $W(S)$,  and is indeed  of zero dimension according to a result of Rudolph. \\
  2.~relies on a classical criterion of Davenport-Erd\" os-Le V\^eque \cite[Lemma 1.8,  p.6]{Bug 2}. 
  
  \begin{prop}\label{del} Let $(\Omega, \mathcal{A}, \mu)$ be a probability space, $(X_n)$ a sequence of complex values random variables with $|X_n|\leq 1$ and let $A_n=\frac{1}{n}\sum_{j=1}^n X_j$ be their averages. Then
 $$\sum_{n=1}^\infty \frac{1}{n} \int_{\Omega} |A_n|^2 d\mu<\infty \Longrightarrow A_n\to 0\quad  \mu \hbox{\ -almost everywhere}.$$
   \end{prop} 
   We will now show, with $X_{n}(x)=e(s_n x)$ that if $\widehat{\mu}(h)=O\big(1/(\log\log h)^\alpha\big)$ for some $\alpha>1$, then $\mu(W(S))=0$. So that $\mu$ can certainly not be supported on $W(S)$! For that, we check the assumptions of Proposition \ref{del}. Clearly, by expansion\\
   $$ \frac{1}{n} \int_{\T} |A_n|^2 d\mu=\frac{1}{n^3} \big(n+\sum_{ j\neq k, 1\leq j,k\leq n} \widehat{\mu}(s_k-s_j)\big)$$
   so that
  $$ \frac{1}{n} \int_{\T} |A_n|^2 d\mu\leq\frac{1}{n^3} \big(n+2\sum_{ 1\leq j<k\leq n} |\widehat{\mu}(s_k-s_j)|\big).$$
  But we know that, for $j<k$, we have
  $s_k-s_j\geq s_k -s_{k-1} \gg s_k/\log s_k$
  so that 
 $$ \log\log (s_k-s_j)\gg \log\log s_k\gg \log (\sqrt {k})\gg \log k$$
 and hence
  $$ |\widehat{\mu}(s_k-s_j)|\ll 1/ (\log k)^\alpha.$$
  This gives us 
  $$ \frac{1}{n} \int_{\T} |A_n|^2 d\mu\ll \frac{1}{n^3} \big(n+\sum_{ 1\leq j<k\leq n}1/ (\log k)^\alpha\big)\ll n^{-2}+ n^{-2} \sum_{k=2}^n 1/ (\log k)^\alpha$$ implying
  $$ \frac{1}{n} \int_{\T} |A_n|^2 d\mu \ll \frac{1}{n(\log n)^\alpha}.$$
 So, the assumption of  Proposition \ref{del} is met, meaning that $\mu$-almost every $x$ satisfies $A_{n}(x)\to 0$. If $q$ is a non-zero integer, we can do the same with   the averages $A_{n}^{(q)}(x)=\frac{1}{n}\sum_{j=1}^n e(qs_j x)$ and conclude that  $\mu$-almost every $x$ is uniformly distributed, i.e. $\mu(W(S))=0$. This clearly ends the proof of Theorem \ref{upr}. We could similarly show that, for $\mu$ carried by $W(S)$, 
 $$\widehat{\mu}(h)=O\Big(\frac{1}{\log \log h(\log_{3}h)^\alpha}\Big),\ \alpha>1$$ is impossible. 
  \end{proof}
  
\noi{\bf Comments.} The property for a set to support a Rajchman measure says more on the   {\it lack of porosity} of this set than on its size. This is perhaps  an illustration of the uncertainty principle for measures: actually, $W(S)$ is rather big 
as it has positive Hausdorff dimension, but so is the porous triadic Cantor set $K$ with Hausdorff dimension $\log2/\log3$ which however supports no Rajchman measure (\cite{KS} p.~59).
 We provide here a direct proof of the last fact based on the following property of $K$ : The Cantor set $K$ is $\sigma_3$-invariant,   in particular, for every $x\in K$ and $k\ge1$, $3^kx\notin I:=]1/3,2/3[$. 
 Suppose that $K$ supports a non-trivial  $M_0$-measure $\mu$. Without loss of generality we assume that $\widehat{\mu}(0) \not=0$ (otherwise we can consider the measure
 $e^{2\pi i at}\mu$ instead of $\mu$, with $a$ such that $\widehat{\mu}(a)\not=0$). Let now $f$ be the sum of an absolutely summable trigonometric series ($f\in A(\T)$), with its closed support inside $I$ and $\widehat f(0)=1$; we then define $f_k\in A(\T)$ by $f_k(t)=f(3^kt)$: for every $k\ge0$, the spectrum of $f_k$ is contained into $3^k\Z$ and its closed support stays away from $K$. It follows that
\begin{equation}\label{obviouss} \forall n\not=0,\ \lim_{k\to\infty}\widehat f_k(n)=0
\end{equation}
and
\begin{equation}\label{obvious}0=\int f_k\ d\mu=\sum_n \widehat f_k(n)\overline{\widehat\mu(n)}.
\end{equation}
Since $\widehat{\mu}(n) \to 0$ as $|n| \to \infty$, given $\varepsilon>0$, there exists $N$ such that $|\widehat\mu(n)|\leq \varepsilon$ for $|n|> N$, implying
$$\big|\sum_{|n|>N} \widehat f_k(n)\overline{\widehat\mu(n)}|\le \varepsilon\sum_n |\widehat f_k(n)\big|= \varepsilon\sum_n |\widehat f(n)|=:C\eps,$$
whence, from (\ref{obvious}),
$$\big|\sum_{|n|\le N} \widehat f_k(n)\overline{\widehat\mu(n)}\big|=\big|\sum_{|n|> N} \widehat f_k(n)\overline{\widehat\mu(n)}\big|\le C\varepsilon.$$
But, $N=N_\eps$ being fixed, (\ref{obviouss}) leads to the contradiction (for $\varepsilon$ small)  
$$0< |{\widehat f(0)}{\overline{\widehat\mu(0)}}|=\lim_{k\to\infty}\big|\sum_{|n|\le N} |\widehat f_k(n)\overline{\widehat\mu(n)}\big| \leq C\varepsilon.$$ 
 
 Another striking example is the set of Liouville numbers which supports  a Rajchman measure (\cite{Blu,Bug}), though of zero Hausdorff dimension.

\subsection{Khinchin class of the Furstenberg set $S$}
From the classical result of the equidistribution mod 1 of $(nx)$, $x\notin\Q$, it follows that ${1\over N}\sum_{n<N}f(nx)\to\int_\T f\ dm$ almost-everywhere for every continuous, or even Riemann-integrable function $f$. Khinchin \cite{Kh} conjectured that the result still holds with any function in $L^1(\T)$. But, Marstrand \cite{Ma} proved that it fails for some function in $L^\infty(\T)$. However the class of $L^1$-functions satisfying Khinchin's conjecture deserves to be explored and, many contributions took this path (see \cite{BW}). Let us cite one of earlier result due to Koksma \cite{Ko}: any $L^2$-function
such that $\sum_{n=1}^\infty |\widehat{f}(n)|^2 (\log \log n)^3<\infty$ satisfies Khinchin's conjecture. 
\medskip

 We are interested in subsequences $(k_n)$ of the sequence of positive integers, especially the Furstenberg sequence $(s_n)$. We propose the following definition.
\begin{defi} The Khinchin class of an increasing sequence $(k_n)_{n\geq 1}$ of integers is the class of functions $f\in L^1(\T)$ satisfying 
\begin{equation}\label{eq:Khintchin}
\lim_{N\to\infty}{1\over N}\sum_{n=1}^N f(k_n x)=\int_\T f\ dm \hbox{\ almost everywhere}.
\end{equation}
\end{defi}
The question of Khinchin, related to other sequences of integers, has been fruitfully studied by Marstrand in the above cited paper \cite{Ma}. He proved the following (among others).
\begin{theo}[Marstrand] \label{jenaimarst}  Every function $f\in L^\infty(\T)$ is in the Khinchin class of the Furstenberg sequence $S$.
\end{theo}
 The assumption $f\in L^\infty(\T)$ was later on dropped by Nair \cite{Na} in a difficult paper: assuming $f\in L^1(\T)$ is enough.  The Khinchin class of $(s_n)$ is thus the whole of  $L^1(\T)$.
 
Using an ergodic argument, we give below a simple proof of Marstrand's result, which nearly recaptures Nair's generalization and holds in the more abstract context of a $\times2\times3$-invariant probability measure, and $S$-ergodic  (if any such measure different from $m$ exists). 
A measurable function $f:\T\to \C$ is in the class $L\log^{+}L$ with $L:=L^1(\mu)$ if 
$$\int_{\T} |f|\log^{+} |f|d\mu<\infty.$$
\begin{theo}\label{nair} Let $\mu$ be a $S$-invariant and $S$-ergodic probability measure. If $f\in L\log^+L$, then, 
$$A_Nf(x):={1\over N}\sum_{n\le N}f(s_nx)\to \int fd\mu\quad  \mu-a.e.$$
In particular
$${1\over N}\sum_{n\le N}\delta_{s_nx}\to \mu \ \ {\rm weakly}^*\quad  \mu-a.e.$$
Consequently, if Furstenberg's conjecture were not correct, there would exist an $S$-invariant  probability measure carried by $W(S)$.
\end{theo}
Observe that $S$-ergodicity means this: $h\circ \sigma_2=h$ \textit{and}  $h\circ \sigma_3=h$ imply $h$ constant $\mu$-a.e.; and this does not infer 2- or 3- ergodicity separately.

\begin{proof} 
 We denote by $\sigma_q$ the $q$-shift: $x\mapsto qx$ on $\T$ and put $\sigma_qf(x)=f(qx)$. Now, assume that   $f\in L^1(\mu)$ and $\int fd\mu=0$; 
then we set 
$M= {\log s_N}/{\log2}$ and we decompose, as in the proof of Proposition \ref{S_N},
 $$\sum_{n\le N}f(s_nx)=\sum_{2^j3^k\le s_N}f(2^j3^kx)=\sum_{j+k\alpha^{-1}\le M}f(2^j3^kx)$$
$$=\sum_{j\le M}\sum_{k\le(M-j)\alpha}f(2^j3^kx)=\sum_{j\le M}\sigma_2^j\big(\sum_{k\le(M-j)\alpha}f(3^kx)\big).$$
Observing that $M\approx \sqrt N$, we get

\begin{eqnarray*}|A_Nf(x)| & \lesssim & \Big|{1\over M^2}\sum_{j\le M}\sigma_2^j\big(\sum_{k\le(M-j)\alpha}f(3^kx)\ \big)\Big|\\
                                           & \lesssim &  \Big|{1\over M}\sum_{j\le M}\sigma_2^j \Big({M-j\over M}F_{M-j}(x)\Big)\Big|
                                         \end{eqnarray*}
where
$$F_K(x):={1\over K\alpha}\sum_{k\le K\alpha}f(3^kx).$$

By Birkhoff's theorem applied to the dynamical  system $(\T,\sigma_3,\mu)$, we know that 
$F_K(x) \to \ell(x)$ $\mu$-a.e. where $\ell=\E_3(f)$ is the conditional expectation with respect to  the 3-invariant sets. 
We thus decompose
$$T_M={1\over M}\sum_{j\le M}\sigma_2^j \Big({M-j\over M}(F_{M-j}(x)-\ell(x))\Big)+
{1\over M}\sum_{j\le M}{M-j\over M}\sigma_2^j \ell(x)
=:T'_M+T''_M.
$$
It is easy to see that the second mean $T''_M \to 0$   by the following remark.
\begin{lem} Let $(z_j)$ be a sequence of complex numbers such that $Z_j/j\to 0$ with $Z_j=z_1+\cdots+z_j$. Then
$$\lim_{M\to\infty}{1\over M}\sum_{j=1}^M\big(1-{j\over M}\big)z_j= 0.$$
\end{lem}
The proof of this lemma is simple:
$${1\over M}\sum_{j=1}^M\big(1-{j\over M}\big)z_j={1\over M^2}\sum_{j=1}^{M-1}Z_j=o\big(\sum_{j=1}^{M-1}j\big)/M^2=o(1). $$
Indeed, we just apply this lemma with $z_j=\sigma_2^j \ell(x)$:  Birkhoff's theorem gives $Z_j/j\to \E_2(\ell)=\E_2\E_3(f)$ $\mu-a.e.$ since $\mu$ is 2-invariant. As $\E_2\E_3(f)$ is $S$-invariant measurable function (observe that both projections $\mathbb{E}_2$ and $\mathbb{E}_3$ are commuting) and $\mu$ is assumed $S$-ergodic,  $\E_2\E_3(f)$ must be constant $\mu$-a.e; so we get $\E_2\E_3(f)=\int fd\mu=0$, then $Z_j/j\to 0$  and finally $T''_M\to 0$.\\

 Now, look at the first mean $T'_M$, which satisfies
$$|T'_M|\leq {1\over M}\sum_{j\le M}\sigma_2^j \big|F_{M-j}(x)-\ell(x)\big|$$
 We will  be done if we can prove that the RHS tends to zero $\mu$-a.e.:
\begin{equation}\label{LlogL}\lim_{M\to\infty}{1\over M}\sum_{j\le M}\sigma_2^j \big |F_{M-j}(x)-l(x)\big|=0\ \ \mu-a.e.
\end{equation}
We will prove (\ref{LlogL}) by making use of the following improvement of Birkhoff's theorem (\cite{Pe}, Ch.VI,  p. 262).
This improvement in \cite{Pe} was not explicitly stated and was proved under the assumption of ergodicity. But this assumption is not necessary.
For clarity  we give a complete proof.  

\begin{lem} Let $T$ be  measure preserving transformation on $(X,\mathcal{B},\mu)$ and let $(G_n)\subset L^1(\mu)$. Suppose that
$G_n\to G$ a.e. and $G^*:=\sup_n|G_n|\in L^1(\mu)$. Then
\begin{equation}\label{eq:SBirkhoff}
     \lim_{n\to \infty} \frac{1}{n}\sum_{k=0}^{n-1} |G_{n-k}- G|\circ T^k =  0   \ \ \mu-a.e
\end{equation}
In particular, we have
\begin{equation}\label{eq:WBirkhoff}
 \lim_{n\to \infty} {1\over n}\sum_{k=0}^{n-1} G_{n-k}\circ T^k = \mathbb{E}(G|\mathcal{I})\ \ \mu-a.e.
 \end{equation}
where $\mathcal{I}$ is the $\sigma$-field of $T$-invariant sets.
\end{lem}

To prove (\ref{eq:SBirkhoff}), we can assume $G=0$ and $G_k\ge0$ for all $k$. Let $g_n= \sup_{k\ge n}G_k$ for $n\ge 1$; by hypothesis, $g_n\le g_1= G^*\in L^1(\mu)$ for every $n$.
Fix $1\le m<N$. We decompose

$$ \frac{1}{N}\sum_{k=0}^{N-1} G_{N-k}\circ T^k  = \frac{1}{N}\sum_{k=0}^{N-m}G_{N-k}\circ T^k  + \frac{1}{N}\sum_{k=N-m+1}^{N-1} G_{N-k}\circ T^k 
$$
$$\le  \frac{1}{N}\sum_{k=0}^{N-m} g_m \circ T^k + \frac{1}{N}\sum_{k=N-m+1}^{N-1} G_{N-k}\circ T^k=:A_N'+A_N''. 
$$

 The term $A''_N$ involves only $m-1$ functions and  satisfies
 $$0\leq \limsup_{N\to \infty} A''_N=  \limsup_{N\to \infty} \frac{1}{N}\sum_{j=1}^{m-1} G_{j}\circ T^{N-j}=0\ \ \mu-a.e.
$$
 (by Borel-Cantelli's lemma, $ \frac{1}{N}\,(h\circ T^N)\to 0\  \mu$-a.e. for every $h\in L^{1}(\mu)$).
Switching to the first term $A'_N$, 
we apply  Birkhoff's theorem to  $g_m$ and  we get
$$
 \limsup_{N\to \infty} A'_N\le    \lim_{N\to\infty}   \frac{1}{N}\sum_{k=0}^{N-m} g_m \circ T^k  = \mathbb{E}(g_m|\mathcal{I}) \ \ \ \ \mu-a.e.
$$
with $\mathcal{ I}$ standing for the $\sigma$-field of $T$-invariant sets.
Finally, we obtain 
$$
   0\le L:= \limsup_{N\to \infty} \frac{1}{N}\sum_{k=0}^{N-1}G_{N-k}\circ T^k\le \mathbb{E}(g_m|\mathcal{I}) \ \ \ \mu- a.e.
$$
 hence
$$0\le \E(L)\le \E(g_m)\ \ \hbox{for every}\ m.$$
  Remember that  $g_m\to 0$ $\mu$-a.s. and $g_m\le G^*\in L^1(\mu)$,  so that the dominated convergence theorem applies: we get $\E(g_m)\to 0$, 
$ \E(L)= 0$
and  $L=0$ $\mu$-pp.
We have thus proved (\ref{eq:SBirkhoff}), which, together with Birkhoff's theorem, implies (\ref{eq:WBirkhoff}):
$$
  \lim_{n\to \infty} \frac{1}{n}\sum_{k=0}^{n-1} G_{n-k}\circ T^k =    \lim_{n\to \infty} \frac{1}{n}\sum_{k=0}^{n-1} G \circ T^k = \mathbb{E}(G|\mathcal{I}) \quad \mu-a.e.
$$ 
Applied with  $T=\sigma_2$ and  $G_n=F_n$, $G=l$, the lemma  gives the result. 
It remains to check that $\sup_K|F_K|\in L^1$ under our assumptions on $f$. And this is nothing but the classical maximal ergodic theorem, valid for any dynamical system $(X,\mathcal{B},\mu,T)$, that we recall below (see \cite{EW}, Ch.II, p. 38, and \cite{Ne2}, Ch.IV,  p. 70-71 for the consequence).
\begin{lem} Given any dynamical system $(X,\mathcal{B},\mu, T)$, $g\in L^1(\mu)$ non-negative and $G_n={1\over n}\sum_{k<n}T^kg$, 
the maximal functions $G_N^\ast=\sup_{n\le N} G_n$ ($N\ge 1$) satisfy
 $$\mu(G_N^\ast>t)\leq \frac{1}{t} \int_{X} g\,{\bf 1}_{\{G_N^\ast>t\}}\;d\mu \hbox{\quad for all}\ t>0.$$
As a consequence, the maximal function $G^\ast=\sup_{n} G_n$
 is in $L^1(\mu)$ as soon as $g$ is in  $L\log^+L$.
\end{lem}
\noi The equality (\ref{LlogL}) is thus proved. Finally, if an  $S$-invariant  probability measure $\mu$ exists with $\mu\neq m$, since $\mu$ is a barycenter of   $S$-invariant  and $S$-ergodic probability measures, we can assume it to be  $S$-ergodic as well, and there exists an integer $k\neq 0$ such that $\widehat{\mu}(k)\neq 0$. Applying  Theorem \ref{nair} to $f=e_{-k}$, we get that, for $\mu$ almost all $x$: 
$$ \frac{1}{N}\sum_{n=1}^{N} e(-ks_nx)\to\int_{\T} e_{-k}d\mu= \widehat{\mu}(k)\neq 0$$
so that $\mu(W(S))=1$. 
\end{proof}

\section{Thin sets in harmonic analysis}
 We first make a short reminder on Orlicz functions and  spaces encountered just before, to be used later in this work. \\
{1.} An Orlicz function $\varphi:[0,\infty[\to \R+$  is a non-negative, increasing and convex function with $\varphi(0)=0$. \\
{2.} If $(X,\mathcal{B},\mu)$ is a probability space and $\varphi$ an Orlicz function, the associated Orlicz space $L^\varphi$ is defined as 
$$L^\varphi=\Big\{f  \hbox{\  measurable}: \int_{X} \varphi\Big(\frac{|f(x)|}{a}\Big)d\mu(x)<\infty\ \hbox{\ for some}\ a>0\Big\}.$$
The Orlicz (Luxemburg)  norm $\|f\|_{L^\varphi}$ of $f$ is the infimum of those $a>0$ such that  $\int_{X} \varphi\big({|f(x)|/a}\big)d\mu(x)\leq 1$. \\
{3.} If $\varphi$ moreover satisfies the $\Delta_2$-condition,  namely $\varphi(2x)\leq C\varphi(x)\  \forall x>0$, then 
$$f\in L^\varphi\Longleftrightarrow \int_{X} \varphi\big(|f(x)|\big)d\mu(x)<\infty.$$
 \noi {\bf Examples.} a) With $\varphi(x)=x^p$ ($1\leq p<\infty$), we recover
 the usual   Banach space $L^p$. \\
b) Take  $\varphi(x)=\varphi_{\alpha}(x):=x\,\log^{\alpha} (1+x),\ 0<\alpha\leq 1$. This function is  an Orlicz function satisfying  $\Delta_2$. Clearly $L^{\varphi_\alpha}\hookrightarrow L^1$.

\smallskip
We next coin combinatorial tools  which reveal adapted to spectral harmonic analysis \cite{LQ}.
\subsection{Dissociate and quasi-independent sets}
\begin{defi} A set $E\subset\Z\backslash\{0\}$ is said to be \textnormal{quasi-independent} if for all \textnormal{distinct} elements $x_1,\ldots,x_n\in E$ and for all $\eps_1,\ldots,\eps_n\in\{-1,0,1\}$, the \textnormal{relation} $\sum\eps_kx_k=0$ implies that all $\eps_k=0$. The quantity $\sum_{k}|\varepsilon_k|$ is called the \textnormal{length} of the relation. 

$E$ is said to be dissociate if for all distinct elements $x_1,\ldots,x_n\in E$ and $\eps_1,\ldots,\eps_n\in\{-2,-1,0,1,2\}$, $\sum\eps_kx_k=0$ implies all $\eps_k=0$.
\end{defi}
A more comprehensible definition is the following:  $E$ is quasi-independent if every finite sum of elements of $E$ admits a unique such decomposition; 
it is dissociate if every expansion $\sum\eps_kx_k$ with $\eps_k\in\{-1,0,1\}$ is unique.

A Hadamard set with ratio $q$ is quasi-independent if $q\ge2$, dissociate if $q\ge3$. Every Hadamard set with ratio $q>1$ is a finite union of dissociate sets.

\subsection{$p$-Sidon sets, $\Lambda(p)$-sets}

 1-periodic functions with spectrum in some Hadamard subset of $\N$ enjoy specific convergence and regularity properties. The most famous ones, that we recall here, have led to define the so-called ``thin sets" in harmonic analysis.  

{\it 1. A continuous function with spectrum in some Hadamard subset $E$ possesses an absolutely convergent Fourier series.}
In other terms, $f\in C_E\Longrightarrow \widehat{f}\in \ell^1$. \ 
The sets of integers enjoying this property have been called {\bf Sidon sets} by Rudin \cite{Ru2} and more generally, for $1\le p<2$,   
$$E\ \hbox{ is a {\bf p-Sidon set}}\Longleftrightarrow [f\in C_E\Longrightarrow \widehat{f}\in \ell^p].
$$
(Of course, 1-Sidon=Sidon.)
If $E$ is $p$-Sidon, there exists a best constant $S_p(E)$, called the $p$-Sidon constant of $E$, such that
$$\Vert \widehat{f}\Vert_{p}\leq S_p(E)\Vert f\Vert_\infty \hbox{\ for every}\ f\in C_E.$$

{\it 2. If $f\in L^1$ has a Hadamard spectrum, then $f$ belongs to all the $L^p$, $p<\infty$.}
The sets $E$ enjoying this property for one $p>2$ have been called {\bf $\Lambda(p)$-sets} (again by Rudin) \cite{Ru2}, and we denote by $\lambda_{p}(E)$ the best constant such that
$$\Vert f\Vert_p\leq \lambda_{p}(E) \Vert f\Vert_2$$ for every polynomial with spectrum in $E$.   Then, every  $f\in L^2_E$ belongs to $ L^p$ and the above inequality still holds. 
See also \cite{LQ} vol.2, p.144-151. 

\medskip
\noi
{\bf First examples.} Quasi-independent sets are Sidon sets \cite{Ru2} with Sidon constant bounded by 8. 
A Hadamard set $E=\{n_j\}\subset\N$ is a $\Lambda(p)$-set for all finite $p$. 
More generally, a Sidon set $E$ is a $\Lambda(p)$-set for all $p>2$ with constant $ \lambda_{p}(E)\leq 2S(E)\sqrt p=O(\sqrt{p})$ (also a result of Rudin, the converse being due to Pisier).
Later on, Bourgain showed that  true $\Lambda(p)$-sets do exist if $p>2$, that means  sets which are $\Lambda(p)$, but not $\Lambda(q)$ as soon as $q>p$. Whence the restriction to the parameters $p>2$.
The question whether true $\Lambda(2)$-sets do exist is an open problem. 
\medskip

We present below characterizations of $p$-Sidon set and $\Lambda(p)$-sets.

\smallskip\noi $\bullet$ There exists a characterization of $p$-Sidon sets in terms of interpolation.
\begin{theo}[\cite{ER}]
Let $1\leq p<2$ and $E\subset\Z$ \textnormal{(}or $\N \textnormal{)}$. Assume $\frac{1}{p}+\frac{1}{q}=1$. The following assertions are equivalent\textnormal{:} \\
\indent {\rm 1)} $E$ is $p$-Sidon.\\
\indent {\rm 2)} $\sum|\widehat f(n)|^p<\infty$ for  $f\in L^\infty_E$.\\
\indent {\rm  3)} If $b\in\ell^q(E)$, there exists $\mu\in M(\T)$ such that $\widehat\mu_{|E}=b$.\\
\indent  {\rm  4)} If $b\in\ell^q(E)$, there exists $f\in L^1(\T)$ such that $\widehat f_{|E}=b$.
\end{theo}

\medskip
\noi{\bf Remark.} We observe that the harmonic classification  is in some sense complementary to the dynamical one detailed in \cite{BL} (see also \cite{Qm} subsection 4.3.2.): the implied set of integers $E$ is viewed as a spectrum in the first classification and as a co-spectrum in the second one. In this respect, we can see that {\it a Sidon set cannot be Hartman-ud}.
Indeed, if $E$ is Hartman-ud, then $\mathcal{C}_E$, the space of continuous functions with spectrum inside $E$, contains a copy of the Banach space $c_0$ of ultimately null sequences (\cite{PIQ}); if $E$ is Sidon, $\mathcal{C}_{E}$ is by definition isomorphic to $\ell^1$, which does not contain $c_0$, hence $E$  will never be Hartman uniformly distributed.

\medskip

\noi$\bullet$ The question whether every Sidon set is a finite union of quasi-independent sets is an open question, motivated by a characterization of a Sidon set due to Pisier (\cite{Pi}): {\it $E\subset\Z\backslash\{0\}$ is a Sidon set if and only if there exists $\delta>0$ such that, from every finite subset $A\subset E$,  a quasi-independent set $B$ can be extracted with $|B|\ge\delta|A|$.} An analogue for $p$-Sidon sets will be considered later (section 5.4).

\medskip

\noi$\bullet$ An arithmetical characterization of these ``thin sets" -- in terms of arithmetical progression or more generally mesh -- is not yet complete (if it ever exists, cf. \cite[Vol.~2,  Ann.~D, p.~316-323]{LQ}), but we have at our disposal necessary conditions, as this one, called {\it mesh condition.}
\begin{defi} For $E\subset\Z$, we denote by $\alpha_E(N)$ the maximal number of elements of $E\cap\{a+b,a+2b,\dots,a+Nb\}$ over all arithmetical progressions  of length $N$.
\end{defi}
\begin{prop}[\cite{Ru2}] \label{maille}
The following holds.\\
\indent {\rm 1)}  If $E$ is $\Lambda(p)$ for $p>2$, then $\alpha_E(N)\le 4 \lambda_{p}(E)^2 N^{2/p}$.

\indent {\rm 2)} If $E$ is $p$-Sidon, then $\alpha_E(N)\le c(\log N)^{p\over 2-p}$.

\end{prop}
\noi {\bf Remarks.} 1. The $\Lambda(q)$-sets with $2\leq q<\infty$ can apparently be less lacunary than $p$-Sidon sets, $1\leq p<2$.  We refer to \cite{Ru2} for the proof of the seminal result {1)}.  \\
2. For $p$-Sidon sets, we have nothing better than for $p$-Rider sets (see Proposition \ref{mc}).
When $p=1$ i.e.~$E$ is Sidon, we can deduce from the result {1)} the estimate of $\alpha_E(N)$ announced in {2)}. In fact,  use  the estimate $\lambda_{p}(E)=O(\sqrt p)$ holding for a Sidon set $E$ to 
 get from {1)} $\alpha_E(N)\le c{p} N^{2/p}$ for every $p>2$.  Let us now choose $p=2\log N$ with $N\ge 3$ so that $p>2$, to get $\alpha_E(N)\le c\log N$ as claimed in {2)}.

\medskip\noi
{\bf An example of cartesian product.} A classical example \cite{Wo}  is $E=\{(\lambda_m, \mu_n)\}=:A\times B$, the cartesian product of two Sidon sets  in $\Z$.  
The set $E$ is exactly $4/3$-Sidon  in the  product group $\Z^2$ \cite{Wo}.
 In particular, $E=\{(2^m, 3^n)\}$ is $4/3$-Sidon in $\Z^2$. 
But the Furstenberg set is a Minkowski product.  As we will see,  it is a $\Lambda(p)$-set (Gundy-Varopoulos), but  we don't know if it is a $p$-Sidon set for some $p$ such that $4/3\leq p<2$.

\medskip
We now present the notion of Paley set, which is related to that  of  $\Lambda(p)$-set. The Furstenberg set is such a set (Gundy-Varopoulos), as we will see.

\subsection{Paley sets} \label{sect:Paley}
 The general framework consists of  a set $E$ of integers and a Banach space $X$ of integrable functions on the circle. We say that \textit{the set $E$ is $X$-Paley if the Fourier transform of any function $f\in X$, once restricted to $E$, is square-summable} i.e.
$$\widehat{f}_{|E}\in \ell^2.$$
In this case, there exists a smallest constant $P(X,E)>0$, {\it the Paley constant of the pair $(X,E)$},  such that the Paley inequality holds:
$$ \Vert \widehat{f}_{|E}\Vert_2\leq P(X,E)\Vert f\Vert_X, \quad \forall f\in X.$$
The terminology comes from a theorem of Paley, concerning specifically the case
$X=H^1,\ E=\{2^n, n\geq 0\}$
 (recall that $H^1=\{ f\in L^1(\T)\ ;\ \widehat{f}(n)=0 \hbox{\ when}\  n<0\}$).
It asserts that the powers of 2 constitute a $H^1$-Paley set.
\begin{theo}[Paley] There exists $C>0$ such that, for every $f\in H^1$,
$$\Big(\sum_{k\ge0}|\hat f(2^k)|^2\Big)^{1/2} \le C\|f\|_{H^1}.$$
\end{theo}

\smallskip
\noindent {\bf Remarks.}
1.~The example of Paley does not extend to $X=L^{1}(\T)$. If $K_N$ denotes the Fej\'er kernel of order $N$, we have 
$$\sum_{n\geq 0}|\widehat{K_N}(2^n)|^2\gtrsim \log N \hbox{\quad but}\quad \Vert K_N\Vert_1=1.$$

\noindent2.~The Paley inequality  holds for a space hardly smaller than $L^1(\T)$, namely the Orlicz space $X=L^{\varphi_{1/2}}$ with our previous notations, and for any Sidon set $E$ \cite[ vol. II, Ch.~XII, p.~132]{Zy}. Hence, it holds for $X=L^q(\T)$ ($q>1$) and Hadamard sets $E$.

\smallskip
\noindent3.~The initial Paley inequality indeed tells that a Hadamard set is $H^1$-Paley; so is any finite union of Hadamard sets.\textit{ But this does not extend to arbitrary Sidon sets}; indeed, Sidon sets $E$ exist for which
$\sup_{N} \big|E\cap [N,2N]\big|=\infty,$
preventing them from being $H^1$-Paley (\cite{LQ} Vol.2, Ch.VI, p.~194-196).

\smallskip 
We will be specially interested in the case $X=L^q(\T)=:L^q$ with $1<q<2$ and will speak of a {\it $q$-Paley set} in that context.
An interesting duality  exists between $q$-Paley sets and $\Lambda(p)$-sets, where $p>2$ is the conjugate exponent of $q$.  
\begin{theo}\label{dual} Let $E$ be a $q$-Paley set for some $1<q<2$. Then $E$ is a $\Lambda(p)$-set and moreover $\lambda_{p}(E)\leq P(L^q, E)$. 
\end{theo}
\begin{proof}  Let $f\in L^{2}_{E}$, and $g\in L^q$ with $\Vert g\Vert_q=1$. We have
$$\Big|\int_{\T} f(x)g(-x)dm(x)\Big|=\Big|\sum_{n\in \Z} \widehat{f}(n)\widehat{g}(n)\Big|=\Big|\sum_{n\in E} \widehat{f}(n)\widehat{g}(n)\Big|$$
$$\leq 
\Vert f\Vert_2\,\big(\sum_{n\in E} | \widehat{g}(n)|^2\big)^{1/2} 
\leq  P(L^q, E)\Vert f\Vert_2 \Vert g\Vert_q= P(L^q, E)\Vert f\Vert_2.$$
 By taking the supremum over the $g$ in the unit ball of $L^q$, we get
\begin{equation}\label{postpal} \Vert f\Vert_p  \leq  P(L^q, E) \Vert f\Vert_2,\end{equation}
hence the result.
 \end{proof}

\smallskip
\noindent {\bf Remark.}  Note the gain when passing from ``$\Lambda(p)$ to $q$-Paley": we enlarge the collection of $f$ involved in the inequality since now the spectrum of the test function $f$ is arbitrary. 

\medskip
Aiming to establish some  properties of thin sets  for the Furstenberg set $S$, we need a detour through probabilities, in particular through martingales.

\section{From Probability to harmonic analysis}
In this paper,  given a probability space $(\Omega, \mathcal{B}, \mathbb{P})$, we will need both \textit{increasing} and  \textit{decreasing}  sequences of sub-$\sigma$-fields $(\mathcal{F}_{n})$ of $\mathcal{B}$. For an expert, there is no real difference (see  \cite[ Ch.IV,  p.124]{Ne1} for a unified presentation), but we prefer to be more specific here. 
\subsection{Azuma's lemmas and subgaussianity of martingales}
 We consider martingales $(M_n)_{0\le n\le N}$ adapted to a finite sequence of sub-$\sigma$-fields $(\mathcal{F}_{n})_{0\le n\le N}$ of $\mathcal{B}$.
$\mathbb{E}_n$ denotes the  conditional expectation given $\mathcal{F}_{n}$, so that, in the increasing case,   
$$\mathbb{E}_{n}(M_{n+1})=M_n,\ 0\leq n<N.$$
We then write $M_N:=\sum_{n=1}^N d_n$, where  
$d_n=M_{n}-M_{n-1},\ \mathbb{E}_{n}(d_{n+1})=0$, and we assume $M_0=0$.
In the case of a decreasing sequence $(\mathcal{F}_{n})$, we will have
$$\mathbb{E}_{n+1}(M_{n})=M_{n+1},\ 0\leq n<N, \ d_n=M_{n}-M_{n+1},\ \mathbb{E}_{n+1}(d_{n})=0.$$ 
A seminal result is:
\begin{theo}[Doob]\label{dodo}  Let $(M_n)_{ n\geq 0}$ be a complex martingale w.r.t. an increasing sequence $(\mathcal{F}_{n})_{n\geq 0}$ of $\sigma$-fields. If $(M_n)$ is bounded in $L^1$, it converges almost everywhere.  
\end{theo}
Suppose now  that $(M_n)$ has bounded increments $(d_n)$. In this context, Azuma proved the following``subgaussian" inequalities \cite[Ch.I, p.~31-32]{LT}.
\begin{lem}[First inequality of Azuma]\label{zuzu}  Let $(M_n)_{0\leq n\leq N}$ be a complex martingale w.r.t.  $(\mathcal{F}_{n})$,  satisfying  $|d_n|\leq C_n$ for $0\leq n\leq N$. When 
$(\mathcal{F}_{n})$ is increasing,   for all $t>0$ we have
$$ \mathbb{P}(|M_N|>t)\leq4\exp\Big(-\frac{t^2}{4\sum_{n=1}^N C_{n}^{2}}\Big)\cdot$$
When 
$(\mathcal{F}_{n})$ is decreasing,   for all $t>0$ we have
$$ \mathbb{P}(|M_0|>t)\leq 4\exp\Big(-\frac{t^2}{4\sum_{n=0}^{N-1} C_{n}^{2}}\Big).$$
\end{lem}
  \noi We will need a more precise version of Lemma \ref{zuzu}, depending on two parameters.
\begin{lem}[Second inequality of Azuma]\label{azuzu}  Let $(M_n)_{0\leq n\leq N}$ be a complex martingale w.r.t. an increasing sequence $(\mathcal{F}_{n})$, which satisfies   $|d_n|\leq a$ and $\sum_{n=1}^N \mathbb{E}_{n-1}(|d_{n}|^{2})\leq b^2$, where $a$ and $b$ are positive constants. 
Then,  for all $t>0$ we have
$$\mathbb{P}\big(\big|M_N\big|>t\big)\leq 4\exp\Big[ -\frac{t^2}{4b^2}\big(2-e^{\frac{at}{b^2}}\big)\Big].$$
In particular, if $\exp(at/b^2)\le3/2$, one has
$$\P(|M_N|>t)\le 4\exp(-{t^2\over8b^2}).$$
Similarly, when 
$(\mathcal{F}_{n})$ is decreasing,   for all $t>0$ we have
$$\mathbb{P}\big(\big|M_0\big|>t\big)\leq 4\exp\Big[ -\frac{t^2}{4b^2}\big(2-e^{\frac{at}{b^2}}\big)\Big].$$
\end{lem}

\subsection{The Salem-Zygmund inequalities}

\noindent The following theorem \cite{SZ} plays a crucial role in Fourier Analysis, and will be deduced here from the first inequality of Azuma. 
We provide short proofs.
\begin{theo}[Salem-Zygmund]\label{anra} Let $(r_n)$ be a Rademacher sequence, $(a_n)$ a scalar sequence, and for all $N\geq 2$, $Q_N$ the random  polynomial 
$$Q_N=\sum_{ n=1}^N r_na_ne_n.$$
Then
$$\E(\|Q_N\|_\infty)\le C\Big(\sum_1^N|a_n|^2\Big)^{1/2}\sqrt{\log N}.$$
\end{theo}
\begin{proof} It holds, using  Bernstein's inequality for trigonometric polynomials of degree $\leq N$ ( \cite{QQ}, Ch.V, p.130; see also \cite{LQ}, vol.1, Ch.VI, p.~241):
$$\E(\|Q_N\|_\infty)\leq C  \mathbb{E}\big(\sup_{x\in E}\vert \sum_{n=1}^N r_n a_n e_{n}(x)\vert\big)=:C\mathbb{E}(\Vert Q_N\Vert_{E}),$$
where $E$ denotes the set of $(4N)-$th roots of unity. The latter supremum is estimated as follows: let $A>1$ and $\sigma_{N}^2=\sum_{n=1}^N |a_n|^2$; then
$$\mathbb{E}(\Vert Q_N\Vert_{E})=\int_{0}^\infty \mathbb{P}(\Vert Q_N\Vert_{E})>t)dt\leq A+\int_{A}^\infty \mathbb{P}(\Vert Q_N\Vert_{E}>t)dt$$
$$\leq A+4N\sup_{x\in E}\int_{A}^\infty \mathbb{P}(\vert Q_{N}(x)\vert>t)dt \leq A+4N\int_{A}^\infty \exp\big(-\frac{t^2}{4\sigma_{N}^2}\big)dt$$
$$\leq A+8N\sigma_N\exp(-\frac{A^2}{4\sigma_{N}^2})$$
(we used Lemma \ref{zuzu}). Now, adjusting $A=\sqrt{4\log N} \sigma_N$
 ends the proof.\end{proof}
 
We will need  a simple fact about the symmetrization \cite{HJ}. Recall that the symmetrization $\tilde X$ of a random variable $X$ is defined (on $\Omega\times \Omega$) by 
$$\tilde{X}(\omega, \omega')=X(\omega)-X(\omega').$$
\begin{prop}\label{raga}  Let $\varphi$ be an Orlicz function, and $X$  a \textnormal{centered} Banach space-valued random variable, with  symmetrization $\tilde X$. Then  we have
$$ \mathbb{E}[\varphi(\|X\|)]\le  \mathbb{E}[\varphi(\|\tilde X\|)]\leq \mathbb{E}[\varphi(2\| X\|)] .$$
\end{prop}

By symmetrization, we get the following more general Salem-Zygmund theorem. 
\begin{theo} \label{sypr}Let $(X_n)_{1\leq n\leq N}$ be a sequence of centered, square-integrable, complex-valued, independent random variables and  $P_N$ the associated random polynomial,  $P_{N, \omega}=\sum_{1}^N X_{n}(\omega)e_n$; then
$$\E(\|P_N\|_\infty)\le C\sqrt{\sum_{1\leq n\leq N} V(X_n)}\sqrt{\log N}.$$
\end{theo}
\begin{proof}Let $\tilde X_n$ be a symmetrization of $X_n$. We write $\sum$ for $\sum_{1\leq n\leq N}$.
 Proposition \ref{raga} gives, for every choice of signs  $(\eps_n)$:
$$\E(\|\sum X_ne_n\|_\infty)\le\E(\|\sum\tilde{X_n}e_n\|_\infty)=\E(\|\sum\varepsilon_n\tilde{X_n}e_n\|_\infty).$$
 We now take $\eps_n=r_n(\omega')$, Rademacher variables independent from the $\tilde{X_n}$, we integrate with respect to $\omega'$ and apply Fubini theorem  to get 
 $$\E(\|\sum X_ne_n\|_\infty)\le\mathbb{E}_{\omega'}\E_{\omega}(\|\sum r_n(\omega')\tilde{X_{n}}(\omega)e_n\|_\infty)$$$$ = \mathbb{E}_{\omega}\E_{\omega'}(\|\sum r_n(\omega')\tilde{X_{n}}(\omega)e_n\|_\infty).$$
Salem-Zygmund's theorem now gives $$\E(\|\sum X_ne_n\|_\infty)\le C \sqrt{\log N}\  \mathbb{E}\ \Big(\sqrt{\sum|\tilde{X_n}(\omega)|^2} \Big),$$
and by  H\"older inequality, 
$$\E\Big(\sqrt{\sum|\tilde{X_n}(\omega)|^2}\Big)\le\Big(\E\sum|\tilde{X_n}(\omega)|^2\Big)^{1/2}=\big(\sum V(\tilde X_n)\big)^{1/2}.$$
Since
$V(\tilde X_n)=2V( X_n)$, we can conclude, suitably changing the constant $C$.\end{proof}

The following specialization will be useful in Section \ref{sect:random}. 

\begin{theo} \label{rec}Let $(\xi_n)_{1\leq n\leq N}$ be  $\{0,1\}$-valued independent variables, with $\mathbb{P}(\xi_n=1)=\delta_n$, and $(\lambda_n)_{1\leq n\leq N}$ scalars of modulus $\leq 1$. Set $m_N=\sum_{n=1}^N \delta_n$. Let  $P_N$ be the random trigonometric polynomial $P_{N,  \omega}=\sum_1^N(\xi_n(\omega)-\delta_n)\lambda_n\,e_n$. Then
$$\mathbb{E}||P_{N}||_\infty\le C\sqrt{m_N\log N}.$$
Moreover, for some constant $C>0$,
$$\mathbb{P}(\Vert P_N\Vert_\infty >t)\leq CN\exp(-\frac{t^2}{8m_N}) \hbox{\ for all}\ t>0 \hbox{\ with}\ e^{t/m_N}\leq 3/2.$$
\end{theo}
\begin{proof} The first part is immediate from Theorem \ref{sypr}, since for the centered  variables $X_n:=\lambda_n(\xi_n-\delta_n)$  we have $V(X_n)\leq\delta_n-\delta_n^2\le\delta_n$. 
For the second part, first use  Lemma \ref{azuzu} with $a=1$ and $b^2=m_N\geq \sum_{n=1}^N |\lambda_n|^2 (\delta_n -\delta_{n}^2)$  to get
$$\mathbb{P}(\vert P_{N}(x)\vert >t)\leq 4\exp(-\frac{t^2}{8m_N}),$$
where $x$ is a fixed $4N$-th root of unity and $t$ as in the statement.
So that, by Bernstein's inequality again,
$$\mathbb{P}(\Vert P_{N}\Vert_\infty >t)\leq CN\exp(-\frac{t^2}{8m_N}).$$
\end{proof}

 \subsection{$p$-Rider sets}
 The functional $\displaystyle{P\mapsto \mathbb{E}(\Vert P_\omega\Vert_\infty)=:[\lbrack P\rbrack]}$ which appeared  in Theorem \ref{anra} (with $P_{\omega}=\sum \varepsilon_{n}(\omega) \widehat{P}(n)e_n$ where $(\varepsilon_n)$ denotes a Rademacher sequence)  is a norm, the ``Pisier norm" on the space $\mathcal{P}$ of all trigonometric polynomials and the completion of  $\mathcal{P}$ with respect to this norm is the so-called Pisier space of almost surely continuous functions (\cite{LQ}, vol.2, Ch.~VI). This norm is much more flexible than the sup-norm because  the characters $e_n$ form an unconditional sequence for this norm,  and motivates the following definition, a variant of $p$-Sidonicity. The set of polynomials with spectrum in a fixed set $E$ is denoted by $\mathcal{P}_E$. 
\begin{defi} A set $E$ of positive integers is called a $p$-Rider set \textnormal{(}with $1\leq p<2$\textnormal{)} if there exists a constant $C$ such that 
$$\Vert \widehat{f}\Vert_p\leq C[[f]] \hbox{\ for all}\ f\in \mathcal{P}_E.$$
The best $C$ is called the $p$-Rider  constant of $E$ and is denoted $s_{p}(E)$. 
\end{defi}
Let us list a few properties of $p$-Rider sets.

\smallskip\noi
\begin{itemize}
\item  A $p$-Sidon set is a  $p$-Rider set, and $s_{p}(E)\leq S_{p}(E)$  (obvious).
\item For $p=1$, we recover the notion of Sidon set, thanks to a result of Rider \cite{Ri}. The analog for $1<p<2$ is an open question. 
\item  If $p<4/3$, a $p$-Rider set is a $q$-Sidon set for each $q>p/(2-p)$ \cite{LR}. 
\end{itemize}
\medskip

As quoted in Proposition \ref{maille}, $p$-Rider sets, in particular $p$-Sidon sets, satisfy an important mesh condition. 
\begin{prop}\label{mc} 
 Let $E$ be a  $p$-Rider set. Then, there exists a constant $C$ such that, for every integer  $N\geq 2$\textnormal{:} 
$$\big|E\cap[1,N]\big|\leq C\big(\log N\big)^{\alpha} \hbox{\ with}\ \alpha=\frac{p}{2-p}\cdot$$
\end{prop}

\begin{proof} Let us write $E\cap[1,N]=\{\lambda_1<\cdots <\lambda_n\}$; set $f=\sum_{j=1}^n e_{\lambda_j}\in C_E$. Using the Salem-Zygmund estimate of Theorem \ref{sypr}, we obtain 
$[[f]]\leq C\sqrt{n}\sqrt{\log N}$. From our $p$-Riderness assumption, this gives
$$n^{1/p}\leq C s_{p}(E)\,n^{1/2} \sqrt{\log N},$$
and the result  ensues. 
\end{proof} 
\noindent {\bf Remark.} Proposition \ref{mc} already shows that the Furstenberg sequence $S$ is at most $4/3$-Rider (in particular  $4/3$-Sidon). Indeed, we saw that  $|S\cap[1,N]|\approx \log^{2}N$, so that, if  $S$ is $p$-Rider, one must have  $2\leq {p}/(2-p)$, or else $p\geq 4/3$. We were not able to show that $S$ is $p$-Sidon or $p$-Rider for $p=4/3$, or even for some $4/3<p<2$. But we shall show that  its randomized version $T$ is $p$-Rider for each $p>4/3$.

  \subsection{Square function of a martingale}
Here we present a fundamental  inequality on the square function of a Hilbert space-valued  martingale due to 
Burkholder  obtained first in 1966 (\cite{Bu}, Theorem 9),
and later optimized in \cite{Bu 2}.
We will use it to prove that the Furstenberg set is a $\Lambda(p)$-set and even a Paley set, results due to Gundy and Varopoulos.  And we will be more precise on the dependence on $p$.

  We consider a martingale $M:=(M_n)_{\ge 0}$ with respect to  a \textit{decreasing} sequence of sub-$\sigma$-fields $(\F_n)_{ n\ge 0}$ of a probability space $(\Omega, \mathcal{F}, \mathbb{P})$.  
  Without loss of generality, we assume that $\F_0=\F$ and $\F_n \downarrow \F_\infty =\{\emptyset, \Omega\}$. Recall that 
\begin{equation}\label{jacque}  
\E_{n+1}(d_n)=0, \quad {\rm with}\ \   d_n:=M_n-M_{n+1}.
\end{equation}
We also assume that $M$ is a Hilbert space-valued  martingale. The square function of the martingale $M$ is the variable defined by 
 $$S=S(M):=\Big(\sum_{n=0}^\infty \Vert d_{n}\Vert^{2}\Big)^{1/2},$$
 where $\|d_n(\omega)\|$ denotes the Hilbert norm.\\
 At this point, we recall a simple  property of $L^p$-spaces (\cite{LiTz} p.46,\cite{LQ} vol. 1 p.189). \\
 $\bullet$ If $p\geq 2$, the Banach lattice $L^p$ is $2$-convex, namely for any $(u_n)\subset L^p$:
 \begin{equation}\label{vexe} \Vert \big(\sum_{n}|u_n|^2\big)^{1/2}\Vert_p\leq \Big(\sum_{n} \Vert u_n\Vert_{p}^{2}\Big)^{1/2}.\end{equation} 
 $\bullet$ If $p\leq 2$, the Banach lattice $L^p$ is $2$-concave, namely
 \begin{equation}\label{cave} \Vert \big(\sum_{n}|u_n|^2\big)^{1/2}\Vert_p\geq  \Big(\sum_{n} \Vert u_n\Vert_{p}^{2}\Big)^{1/2}.\end{equation}

\subsubsection{Burkholder's inequality.}
Let $H$ be a Hilbert space. 
A given $H$-valued integrable  function $f\in L^1(H)$  produces a special martingale $M(f):= (\E_n(f))$. We have the decomposition 
\begin{equation}\label{eq:decomp}
f -\E f =\sum_{n=0}^\infty d_n(f).
\end{equation}
We state Burkholder's inequality for such martingales. 
  In the following, $q$ denotes the conjugate exponent of  $p$ and $p^\ast=\max(p,q)$. It holds (see \cite{Bu 2}, Theorem 3.1 p. 87) 
\begin{prop} [Burkholder inequality]\label{hilb} Let $H$ be a separable  Hilbert space and  let $f\in L^{p}(H)$ with zero mean ($1<p<\infty$). For the martingale  $M_n=\E_{n}(f)$ with its martingales differences  $d_n:=d_n(f)=M_n-M_{n+1}$,\
we have
\begin{equation}
\label{BS} B_p^{-1}\|f\|_p\le\|S(f)\|_p\le B_p\|f\|_p
\end{equation}
where $S(f)^2={\sum_{n=0}^\infty \Vert d_n\Vert^2}$ and $B_p=p^\ast-1$.   
\end{prop}

\subsubsection{The Furstenberg set $S$ is $\Lambda(p)$.}
We will first prove  that the  Furstenberg set $S$ is a $\Lambda(p)$-set, a special case of the  forthcoming  Gundy-Varopoulos theorem.  

The proof will use a special kind of scalar martingales. Let us first recall a basic fact.
Let $m\ge 2$ be an integer. Denote by $\F^{(m)}$ the $\sigma$-field of $m^{-1}$-periodic  Borel sets in $\mathbb{T}$.  
The following basic relations are easy to check.

\begin{lem}\label{lem:basic}  For any $f\in L^1(\mathbb{T})$, we have
\begin{equation}\label{lem:m-cond}
      \mathbb{E}(f| \F^{(m)})(x) = \frac{1}{m}\sum_{k=0}^{m-1} f\Big(x +\frac{k}{m}\Big) = \sum_{m|n} \widehat{f}(n) e_n(x).
\end{equation}
We can also write $ \mathbb{E}(f| \F^{(m)})(x) = f\ast \omega_m (x) $ with  $\omega_m = m^{-1}\sum_{k=0}^{m-1}\delta_{k/m}$.
\end{lem}

Observe that $\F^{(m)} = \sigma_m^{-1}(\mathcal{B}(\mathbb{T}))$ where $\sigma_m$ is the map defined $\sigma_m(x) =m x \mod 1$. The first equality in  (\ref{lem:m-cond})
can be stated as $ \mathbb{E}(f| \sigma_m^{-1}(\mathcal{B}(\mathbb{T})))(x)= R\circ \sigma_m(x)$ where $R$ is the transfer operator
$$
    Rf(x) =  \frac{1}{m}\sum_{k=0}^{m-1} f\Big(\frac{ x + k}{m}\Big)
$$
associated to the measure-preserving dynamical system $(\mathbb{T}, \mathcal{B}(\mathbb{T}), {\rm \lambda}, \sigma_m)$, $\lambda$ being the Haar measure (cf. \cite{Fan2017}, Lemma 2.4).

Now we fix an integer $a\ge 2$. For $n\ge 0$, let $ \F_n =  \F^{(a^n)}$ be the $\sigma$-field of $a^{-n}$-periodic Borel sets in $\mathbb{T}$. From the dynamical point of view, the iteration of $\sigma_a$ is involved. Clearly 
$ \F_0$ is the Borel field $\mathcal{B}(\mathbb{T})$ and $ \F_n$ decreases to the trivial $\sigma$-field $\{\emptyset, \mathbb{T}\}$.
When $g$ is a trigonometric polynomial with zero mean, we can decompose $g$ relatively to the $\sigma$-fields $ \{\F_n\}$
according to (\ref{eq:decomp}). For simplicity, the corresponding martingale difference is still denoted by 
$d_ng$ and the corresponding square function by $Sg$. 
\begin{lem}\label{fiex} Denote by $L_a$ the set of integers not divisible by $a$. Then
\begin{equation}\label{follett} d_{n}g(x)=\sum_{l\in L_a} \widehat{g}(a^{n} l)e(a^{n}lx).\end{equation}
Hence
\begin{equation}\label{ken}S^2(g)(x)=\sum_n\big |\sum_{l\in L_a}\widehat g(a^nl)e(a^nlx)\big|^2.\qquad\qquad \end{equation}
\end{lem}
\begin{proof} By Lemma \ref{lem:basic}, $\E_{n}(g)(x)= g\ast\omega_{a^n}(x)$
 with $\widehat\omega_n={\bf 1}_{a^n\Z}$. Since $d_ng=\E_{n}(g)-\E_{n+1}(g)$,
 (\ref{follett}) follows from
$\widehat {d_ng}=\widehat g\cdot({\bf 1}_{a^n\Z}-{\bf 1}_{a^{n+1}\Z})= \widehat g\cdot{\bf 1}_{a^n L_a}$.
\end{proof}

We denote by $S(a_1,\dots, a_s)$  the semi-group generated by $s$ integers $a_i\ge 2$. 
We say that $a_1,\dots, a_{s}$  are {\it multiplicatively independent}, if
the decomposition of any $d\in S(a_1,\ldots, a_s)$ as $d=a_1^{n_1}\cdots a_{s}^{n_{s}}$ is unique (e.g. $a_1=12,\ a_2=18$). 
Notice that $S(a_1,\ldots, a_s)\subset S(p_1,\ldots, p_t)$ where the $p_j$'s are the primes involved in the decomposition of $a_1,\ldots, a_s$.
It is possible that  $t\neq s$ (e.g. for $a_1=24,\  a_2=30$, we have $s=2, t=3$).

\begin{theo}[Gundy-Varopoulos] \label{bernadette}  Suppose that $a_i \ge 2$ ($1\le i\le s$) are $s$
coprime integers. 
The set   $S(a_1, \dots, a_s)$ is $\Lambda(p)$ for all $p>2$, i.e. $\|f\|_p\le C_{s}(p)\|f\|_2$ if $f\in \mathcal{P}_{S(a_1, \dots,a_s)}$. Moreover, $C_{s}(p)\le C p^{s-1/2}$ with $C$ a constant.

\end{theo}
\begin{proof} We prove the result by induction  on $s\ge1$. The set $S(a_1)= \{a_1^n: n \ge 0\}$ is $\Lambda(p)$ because it is a Sidon set, and then we have
$\lambda_p(S(a_1))= O(\sqrt{p})$.  Suppose that $S(a_1, \dots, a_{s-1})$ is $\Lambda(p)$, namely for $f\in \mathcal{P}_{S(a_1, \dots, a_{s-1})}$ we have
\begin{equation}\label{eq:s-1}
 \|f\|_p \le C_{s-1}(p)\|f\|_2
\end{equation} 
for some constant $C_{s-1}(p)$. 
Assume $f\in\mathcal{P}_{S(a_1, \dots, a_s)}$.  We  distinguish $a_s=a$ and write
$$.$$
\begin{equation}\label{peec} 
f(x):= \sum_{n}  \sum_{m_1, \dots, m_{s-1}} a_{m_1, \dots, m_{s-1}, n} e( a_1^{m_1} \dots a_{s-1}^{m_{s-1}}a^n x)
=\sum_{n} f_{n}(a^n x),
\end{equation}
with
$$
 f_n(y)=   \sum_{m_1, \dots, m_{s-1}} a_{m_1, \dots, m_{s-1}, n} e( a_1^{m_1} \dots a_{s-1}^{m_{s-1}} y) \in \mathcal{P}_{S(a_1, \dots, a_{s-1})}.
$$
Since $a_1, \dots, a_s$ are coprime, no number of $S(a_1, \dots, a_{s-1})$ is divisible by $a$. Thus, 
applying Lemma \ref{fiex} we get $d_{n}(f)(x)=f_n(a^n x)$.
 The  associated square function  is then equal to 
$$Sf(x):=\Big(\sum_n|d_n(x)|^2\Big)^{1/2}=\Big(\sum_n|f_n(a^nx)|^2\Big)^{1/2}.$$
By the scalar Burkholder inequality (\ref{BS}), 
$\Vert f\Vert_{p}\le B_{p}\Vert Sf\Vert_{p}$. It remains to estimate $\Vert Sf\Vert_{p}$.
An important point is that the functions $d_n$ have disjoint Fourier spectra  (because $\hbox{\ sp} (d_n)\subset a^nS(a_1, \dots, a_{s-1})$),  and hence are orthogonal in $L^2$. 
Firstly, the induction hypothesis  (\ref{eq:s-1}) implies
  \begin{equation}\label{I}
      \Vert d_n\Vert_{p}=\Vert f_n\Vert_{p}\le C_{s-1}(p) \Vert f_n\Vert_2= C_{s-1}(p)  \Vert d_n\Vert_2.
\end{equation} 
Then we have by (\ref{vexe}), since $p>2$: 
\begin{eqnarray*}\Vert Sf\Vert_{p}
 & = &\big\Vert\big(\sum |d_n|^2\big)^{1/2}\big\Vert^{}_{p}
 \le 
\Big(\sum \Vert d_n\Vert_{p}^{2}\Big)^{1/2}\\
&\le& C_{s-1}(p) \Big(\sum \Vert d_n\Vert_{2}^{2}\Big)^{1/2}
= C_{s-1}(p) \Big( \big\Vert \sum d_n\big\Vert_{2}^{2}\Big)^{1/2}\\
&=& C_{s-1}(p) \Vert f\Vert_2.
\end{eqnarray*}
We have  used  inequality (\ref{I}), the orthogonality of $d_n$'s and the decomposition $f=\sum d_n$.  Finally,
$$
  \Vert f\Vert_{p}\le B_{p}\Vert Sf\Vert_{p} \le  B_{p}C_{s-1}(p) \Vert f\Vert_{2}.
$$
Thus we have the claimed result with $C_s(p)= B_{p}C_{s-1}(p)=B_p^{s-1} C_1(p) = O(p^{s-1} \sqrt{p})$ by Burkholder's inequality.
This ends the proof of Theorem \ref{bernadette}.
\end{proof}

\subsection{The Gundy-Varopoulos theorem: $S$ is Paley}

Gundy and Varopoulos (\cite{GV}) proved  the following  theorem of  \textit{Paley type} about the semi-group generated by $s$ prime numbers $D_s:=S(p_1, \cdots, p_s)$. 
This  theorem goes beyond the Hadamard case considered by Paley (see Section \ref{sect:Paley}) 
and beyond the previously proved $\Lambda(p)$-character of  $D_s$ (see Theorem \ref{dual}). 
The proof will be essentially the same as for Theorem \ref{bernadette}, except that we will need {\em two} martingales, one of them being vectorial, whereas only one scalar martingale was needed for proving the $\Lambda(p)$-character of $D_s$. 
\begin{theo}[Gundy-Varopoulos]\label{helene}The set $D_s$ is a $q$-Paley set for all $1<q<2$, or else, for every function $f\in L^q $, 
\begin{equation}\label{23}\big(\sum_{n\in D_s} |\widehat{f}(n)|^{2}\big)^{1/2}=\big(\sum_{n_j\in \N_0}|\widehat {f}(p_1^{n_1}\cdots p_s^{n_s})|^2\big)^{1/2}\le  B_{s}(q)\Vert f\Vert_q.
\end{equation}
Moreover $B_{s}(q)\leq  (p-1)^{s}$ where  $p=q/(q-1)$.
\end{theo}

\noi {\em Proof.}
We will first assume that $s=2$ to ease notations. So, we consider two 
distinct primes $a\ge 2$ and $b\geq 2$. We want to show that  the multiplicative semi-group $S(a,b)$ generated by $a$ et $b$ is $q$-Paley. We will denote by   $L:=L_{a,b}\subset S(a,b)^c$ the set of integers $\geq 1$ which are divisible neither by $a$ nor by $b$. 
The proof is based on the following observations.
\\
\noindent  $\triangleright$  We can assume that $f$ is a polynomial  and $\widehat{f}(0)=0$. Our aim is to dominate the quantity  
  \begin{equation}\label{weset}\sum_{m,n}|\widehat{f}(a^m b^n)|^2=:\sigma^2
  \end{equation}  
  by $\|f\|_q^2$ up to a multiplicative constant. Our  first observation is that $\sigma$ is bounded by $\|\mathcal{S}\|_q$
  where $ \S=\big(\sum_{m,n}|g_{m,n}|^2\big)^{1/2}$ with 
$$g_{m,n}(x)=\sum_{l\in L}\widehat f(a^mb^nl)e(a^mb^nlx), \quad (n\ge 0, m\ge 0). $$
Indeed, to let $L^q$-norms appear, we observe that
$$|\widehat{f}(a^m b^n)|^2=|\widehat g_{m,n}(1)|^2\leq \Vert g_{m, n} \Vert _{q}^{2}.$$
Summing over  $m,n$ and using (\ref{cave}) for the exponent $q<2$, we get
\begin{equation}\label{evident} \sigma^2 \leq \sum_{m,n} \big\Vert g_{m,n} \big\Vert _{q}^{2}\leq \Vert \S\Vert_{q}^{2}.
\end{equation} 
We will interpret $\S$ as the square function of a suitably chosen martingale,  which will be handled thanks to Burkholder's inequalities.
\medskip

\noindent $\triangleright$
Let 
$\F^{(m)}$ (resp. $\G^{(n)}$) be the $\sigma$-field of $a^{-m}$-periodic  (resp. $b^{-n}$- periodic) Borel sets, as before.
We write $\E_m=\E(\cdot|\F^{(m)})$ and  $\E'_n=\E(\cdot|\G^{(n)})$ the corresponding conditional expectations (which  commute):
$$\E_m\E'_n=\E'_n \E_m.$$
Indeed,
for $f\in L^1$, one has
$$\E_{m}\E'_{n}f(x)=\E'_{n}\E_{m}f(x) 
=f\ast\omega_{b^{-m}}\ast\sigma_{a^{-n}}(x)=f\ast\sigma_{a^{-n}}\ast\omega_{b^{-m}}(x),$$
by Lemma \ref{lem:basic}.
We accordingly write $d_m=\E_m-\E_{m+1}$ and  $d'_n=\E'_n-\E'_{n+1}$.
 Let $f\in L^1$ and consider its  Fourier expansion
$$f(x)=\sum_{\mu,\nu,\atop l\in L} a_{\mu,\nu, l}\, e(a^{\mu}b^{\nu}lx).$$
As Lemma \ref{fiex} shows, passing from $f$ to $d'_{n}(f)$ consists in freezing $\nu=n$ in this expansion, and then passing to $d_md'_{n}(f)$  in freezing $\mu=m$ as well. So,
the above commutativity  leads to 
\begin{equation}\label{freeze} g_{m,n}=d_m(d'_{n})=d'_{n}(d_m), \qquad \ \mathcal{S}^2=\sum_{m,n} |d_m(d'_n)|^2=\sum_{m,n} |d'_n(d_m)|^2.\end{equation}
(Here and in what follows we omit the letter $f$).
\medskip

 \noindent $\triangleright$ 
  If  $F$ is $H$-valued, we will write $\|F\|_q$ for $\Vert F\|_{L^q(H)}=(\int\|F\|_H^q\,d\mathbb{P})^{1/q}$. \\
  \noindent  Let $H=\ell^2$.  Consider  the vector-valued function  $F\in L^{q}(\ell^2)$ defined by
$$F=(d_0,d_1,\dots,d_m,\dots) : \Omega\to \ell^2.$$
Then consider the vector-valued martingale 
$M_n=\E'_{n}(F)$ with its martingale difference $d'M_n=\E'_{n}(F)-\E'_{n+1}(F)$. 
Observe that
$$d'M_n=(d'_n(d_0),\ d'_n(d_1),\dots, d'_n(d_m),\dots) : \Omega\to \ell^2,$$
and that, by (\ref{freeze}),  the square function associated to $(M_n)$  is equal to
\begin{equation}\label{sq} \Sigma^2:=\sum_{n} \Vert d'M_n\Vert^2=\sum_{m,n} |d'_{n}(d_m)|^2, \hbox{\ or else}\ \Sigma=\mathcal{S}.
\end{equation}
$\triangleright$ On the other hand, the square function of the scalar martingale $(\E_m(f))$, which is  defined by $Sf=\big(\sum_{m}|d_m|^2\big)^{1/2}$, is nothing but
\begin{equation}\label{ao} Sf(\omega) = \Vert F(\omega)\Vert_{\ell^2}.
\end{equation}

  Now, we end the proof as follows, using twice the Burkholder inequalities via (\ref{evident}),  (\ref{sq}) and (\ref{ao}):
$$\sigma\le\Vert \mathcal{S}\Vert_q=\Vert \Sigma\Vert_q\leq B_q\Vert F\Vert_q=B_q\Vert Sf\Vert_q\leq B_{q}^{2} \Vert f\Vert_q\le (p-1)^2\Vert f\Vert_q.$$

 \noindent $\triangleright$ The proof  extends inductively to $s\geq 3$.  To treat the general case of $s$ distinct primes $p_1,\ldots, p_s$, we need some additional notations. 
  If $\alpha:=(n_1,\ldots, n_s)\in \N_{0}^{s}$ and $n=p_{1}^{n_1}\cdots p_{s}^{n_s}$, we write  $n=p^\alpha$. We denote $E_{n}^{(k)}$ the conditional expectation w.r.t.  the $\sigma$-field of $p_{k}^{-n}$-periodic sets, and $d_{n}^{(k)}$ the corresponding increment, as well as 
  $$\delta_{\alpha}=d_{n_s}^{(n_s)}\circ d_{n_{s-1}}^{(n_{s-1})}\circ\cdots\circ d_{1}^{(n_1)}.$$    Now, fix a trigonometric polynomial $f$, write for short $\delta_\alpha=\delta_{\alpha}(f)$; we set 
   $$
   \sigma_s=\big(\sum_{\alpha\in \N_{0}^{s}}|\widehat{f}(p^\alpha)|^2\big)^{1/2},\ {\rm and}\  \Sigma_s=\big(\sum_{\alpha\in \N_{0}^s} |\delta_{\alpha}(f)|^2\big)^{1/2}
   $$  
 the multivariate square function. 
 We prove that $\sigma_s\leq (p-1)^s \Vert f\Vert_q$ in two steps.\\
 {\it Step 1.} We have (a) $\sigma_s\leq  \Vert  \Sigma_s\Vert_q,$
 with the same proof as for $s=2$.\\
  {\it Step 2.} We have (b) $\Vert \Sigma_s\Vert_q\leq (p-1)  \Vert  \Sigma_{s-1}\Vert_q$.
 
  We just need to observe that $\Sigma_s$ is the square function associated to the martingale $(M_{n_s}(F))$, where $F:\Omega\to \ell^{2}( \N_{0}^{s})$ is defined by 
  $F(\omega)=(\delta_{\beta}(\omega))_{\beta}$; here $\beta=(n_{1},\ldots, n_{s-1})$ runs over $\N_{0}^{s-1}$. Then, Burkholder's inequality gives (b).  Clearly, Theorem \ref{helene} follows from (a) and (b).

\medskip

The conclusion of Theorem \ref{helene} holds for $S(a_1,\ldots, a_s)$ when $a_j's$ are multiplicatively independent, because 
$S(a_1,\ldots, a_s)\subset S(p_1,\ldots, p_t)$ where the $p_j$'s are the prime factors of  $a_1,\ldots, a_s$.
One has then to abandon  a precise estimate for the constant $B_{s}(q)$ in terms of $s$ and $q$. 


\section{Random Furstenberg set}\label{sect:random}

We begin by a definition, already coined by Erd\" os and R\'enyi (\cite{Er}), and then systematically studied by Bourgain (\cite{Bo}), which became popular under the name ``Selectors of Bourgain".  We mention in this respect the papers \cite{Ne}, \cite{LQR}, \cite{LQR2}  and \cite{KK1}, \cite{KK2}, \cite{FS}, \cite{BGM}. Let $(\delta_k)_{k\geq 1}$ be a sequence of numbers such that $0<\delta_k<1$ for each $k$, and   $(\xi_k)_{k\geq 1}$ be a sequence of  independent  Bernoulli variables, defined on some  probability space $(\Omega,\P)$, with  
$$\E(\xi_k)=\mathbb{P}(\xi_k=1)=\delta_k,\ k=1,2,\ldots, \hbox{\ and}\  R=R(\omega)=\{k\geq 1\  ;\  \xi_{k}(\omega)=1\}.$$
The set $R$ is called the random set of integers associated with the sequence $(\delta_k)$.
We will always  assume in this section that  
\begin{equation}\label{aa}(\delta_k) \hbox{\ is non-increasing\  and}\  \sum_{k=1}^\infty \delta_k=\infty.\end{equation}
Then, by the Borel-Cantelli lemma, $R$ is almost surely infinite, and we write 
$R=\{u_1<u_2<\cdots<u_n<\cdots\}.$  
 We will  also set
$$m_N=\sum_{k=1}^N \delta_k=\E(\xi_1+\xi_2+\cdots +\xi_N)=\E(|R\cap[1,N]|)=:\E(|R_N|).$$
We will mostly  assume that the ``Bourgain condition'' holds: 
\begin{equation}\label{ma} m_N/\log N\to \infty.
\end{equation}   
 (this is the case when $k\delta_k\to \infty$).
 The question will next be: how to choose $\delta_k$ in order that $R(\omega)$ appears as a reasonable random version of the  Furstenberg set $S$? And then what harmonic analysis and distribution properties does it almost surely possess?

\smallskip 
The  set $S=(s_n)_{n\geq 1}$ satisfies the lacunarity property $|S\cap[1,N]|\approx (\log N)^2$, or again $s_N\approx ae^{b\sqrt N}$. This motivates the choice (to which we will stick)
 \begin{equation}\label{choi} \delta_k=\frac{\log k}{k}.
 \end{equation} 
  In fact, one then has
 \begin{equation}\label{thha} m_N=\sum_{k=1}^N\frac{\log k}{k}= \int_{1}^{N} \frac{\log t}{t}dt+ O(1)=\frac{1}{2} \big(\log N\big)^2+  O(1).\end{equation}
 Notice that $\delta_k=O(1/k)$ corresponds to a random Sidon set (\cite{KK1} p.~364).
We will then reserve the letter $T$ (instead of $R$) for this {\it random version (\ref{choi}) of the Furstenberg set $S$}. 

  The  set $S$ is quite rigid from the arithmetical point of view.  We were partially able to deal  with this arithmetic, and proved  that $S$ is $\Lambda(p)$ for all $p>2$, is not Hartman uniformly distributed, and not $p$-Rider for $p<4/3$.  We will see that the random model $T(\omega)$, our main interest here, almost surely satisfies this $p$-Riderness property for $4/3< p<2$,  among others.

\subsection{$T$ is almost surely weakly lacunary}

We begin by a weak lacunarity property of general random sets $R$ of integers. Then, we switch to our main case of interest: the random Furstenberg set $T$.

\begin{theo}\label{diviav} Let $R=(u_n)$ be the random set of integers corresponding to the sequence $(\delta_k)$ as in (\ref{aa}). Then\\
\indent {\rm 1)}  $R$ satisfies almost surely: $\limsup_{n\to \infty} (u_{n+1}-u_n)=\infty$.\\
\indent {\rm 2)}  Assume moreover that $\sum_{k}\delta_{k}^{2}<\infty$.  Then, $R$ is almost surely weakly lacunary, i.e.  $\lim_{n\to \infty}(u_{n+1}-u_n)=\infty.$    
\end{theo}
\begin{proof} 1)  Fix an integer $q\geq 1$. Consider the events 
$$
A_k=\{ \xi_{kq}=1,\ \xi_{kq+j}=0,\ 1\leq j\leq q-1\},\  \hbox{\ which satisfy}\  
$$
 $$  
 \mathbb{P}(A_k)=\delta_{kq}\prod_{j=1}^{q-1}(1-\delta_{kq+j})\geq\delta_{kq}\prod_{j=1}^{q-1}(1-\delta_{j})=:c_q\delta_{kq}.
 $$
Hence, $\sum_{k\geq 1} \mathbb{P}(A_k)=\infty.$ Since the $A_k$'s are independent, Borel-Cantelli's lemma gives $\mathbb{P}(A)=1$ where
$A=\limsup_{k} A_{k}$.   
Take $\omega\in A$, let $k$ large with  $\xi_{kq}(\omega)=1$ and $n$ the random (large) integer such that $u_n=kq$. Then, $u_{n+1}-u_n\geq q$. This shows that $\limsup (u_{n+1}-u_n)\geq q$ a.s. So that  $\limsup (u_{n+1}-u_n)=\infty$ a.s.

\smallskip
2) Let $(q_k)$ be a sequence of positive integers with $q_k\to \infty$ and $\sum_{k\geq 1}q_k\delta_{k}^{2}$ $<\infty.$
  For a pair $1\leq k<l$ of positive integers, consider the events  
 $$F_{k,l}=\{\xi_k=1,\ \xi_l=1\},\ G_k=\bigcup_{k+1\leq l\leq k+q_k} F_{k,l},\ B=\liminf_{k} G_{k}^{c}.$$
We have $$\mathbb{P}(F_{k,l})=\delta_k \delta_l\leq \delta_{k}^{2},\quad  \mathbb{P}(G_k)\leq q_k\delta_{k}^{2}.$$
 This implies $\mathbb{P}(B)=1$; in other terms, we have almost surely for $k$ large:
$\xi_{l}=0$ whenever $\xi_k=1$ and $k<l\leq k+q_k$. Taking $k=u_n$ gives $u_{n+1}-u_n\geq q_k$ and the result.
\end{proof}

We now specialize to  the random Furstenberg set $T=:\{t_n\}$, corresponding to the choice $\delta_k=\log k/k$, to get  more quantitative results. Notably,  $T$ shares some precise gap properties with its deterministic relative $S$.

\begin{theo}\label{tijdeman} Almost surely, the difference $t_{n+1}-t_n$ satisfies\textnormal{:}\\ 
\indent {\rm 1)}  $\displaystyle\limsup_{n\to \infty} \frac{t_{n+1}-t_n}{(t_n/\log t_n)\,\log\log t_n}\leq 2;$\\
\indent {\rm 2)} $\displaystyle\liminf_{n\to \infty} \frac{t_{n+1}-t_n}{t_n/(\log t_n)^{3+\delta}}\geq 1 \hbox{\quad for all}\ \delta>0.$\\

In particular, the weakly lacunary set $T$ satisfies: almost surely $t_{n+1}/t_n\to 1.$

\end{theo}
\begin{proof} 1) Let $(p_k)$ be a sequence of positive integers such that $p_k=o(k)$. Consider the event 
$E_k=\{ \xi_k=1,\  \xi_{k+1}=\cdots =\xi_{k+p_k}=0\}.$
We see that, given $\varepsilon>0$ and $k$ large (depending on $\varepsilon$) 
$$\mathbb{P}(E_k)=\frac{\log k}{k}\prod_{j=k+1}^{k+p_k} \big(1-\frac{\log j}{j}\big)\leq \frac{\log k}{k}\exp\big(-\sum_{j=k+1}^{k+p_k} \frac{\log j}{j}\big) $$$$
\leq  \frac{\log k}{k}\exp\big(-p_k \frac{\log (k+p_k)}{k+p_k}\big)\leq \frac{\log k}{k} \exp\big(- (1-\varepsilon)p_k\frac{\log k}{k}\big).$$
We now take (ignoring the integer part issues) 
$$p_k=a \frac{k}{\log  k}\log \log k \hbox{\quad with}\ a=\frac{2+\varepsilon}{1-\varepsilon}$$
to get (for $k$ large)
\begin{equation}\label{boca} \mathbb{P}(E_k)\leq   \frac{\log k}{k}\frac{1}{(\log k)^{a(1-\varepsilon)}}=\frac{1}{k(\log k)^{1+\varepsilon}}.\end{equation}
So that $\mathbb{P}\,(\liminf E_{k}^{c})=:\mathbb{P}(A)=1.$
Now, fix $\omega\in A$. For $k\geq k_{0}(\omega)$, one has $\omega\in  E_{k}^{c}$. Take $k=t_n$ with $n$ large enough to ensure $k\geq k_{0}(\omega)$. Then, $\xi_k=1$, and hence $\xi_{j}=1$ for some $k<j\leq k+p_k$, meaning that 
$t_{n+1}\leq k+p_k=t_n +p_k.$
 It easily ensues, for $\omega\in A$:   
 $$\limsup_{n\to \infty} \frac{t_{n+1}-t_n}{(t_n/\log t_n)\,\log\log t_n}\leq  a=\frac{2+\varepsilon}{1-\varepsilon}.$$
 Letting $\varepsilon$ tend to zero, we get the first assertion. \\
 \indent 2) The second one is similar: set 
 $q_k=k/(\log k)^{3+\delta},\hbox{\ where}\  \delta>0$,
 and for  $1\leq k<l$ positive integers, 
$$F_{k,l}=\{\xi_k=1,\ \xi_l=1\},\ F_k=\bigcup_{k+1\leq l\leq k+q_k} F_{k,l}.$$
We get $ \mathbb{P}(F_{k})\leq c^2 q_k \big(\frac{\log k}{k}\big)^2=\frac{c^2}{k(\log k)^{1+\delta}}.$
Once again,  $\mathbb{P}(\liminf_{k} F_{k}^{c})=1$, and almost surely, for $n$ large enough 
$$ t_{n+1}-t_n\geq\frac{t_n}{(\log t_n)^{3+\delta}}\cdot$$
 \end{proof}
\subsection{$R$ is almost surely Hartman uniformly distributed}

The following random result appeared in \cite{Bo} (see also \cite{Ne}), with different terminology.  
 \begin{theo}[Bourgain]\label{BH} Let $(\xi_k)$ be a sequence of  independent random $0-1$-valued variables as before, where $(\delta_k)$  satisfies the Bourgain condition (\ref{ma}).
 Then, a.s. the sequence $R:=\{k\ ;\ \xi_k=1\}$ is Hartman uniformly distributed.
 \end{theo}
\begin{proof}We give a  proof  relying on Section 4 (Salem-Zygmund theorem).
To that effect, we decompose the sums
$$\sum_{n\in R_N}e(nx):=\sum_{n\le N}\xi_ne(nx)=\sum_{n\le N}(\xi_n-\delta_n)e(nx)+\sum_{n\le N}\delta_ne(nx)
$$
and  the deterministic sums $\sum_{n\le N}\delta_ne(nx)$, when $x\in \T,\  x\not=0$, is fixed, are uniformly bounded in $N$ (since $(\delta_n)$ is non-increasing) by a constant $C(x)<\infty$.
As for the first random sum $P_N$, we resort to  Theorem  \ref{rec}.\\
Choose  $t= 6\sqrt{m_N \log N}$. We have $t/m_N=o(1)$ since $m_N/\log N\to \infty$, and thus $\exp({t}/{m_N})\leq 3/2$ for $N$ large enough.  Moreover, $t^2/(8 m_N)\geq 4\log N$. Theorem \ref{rec}  now gives, for $N$ large enough,  
\begin{equation}\label{devi} \mathbb{P}\big(\Vert P_N\Vert_\infty> 6 \sqrt{m_N\log N}\big)\ll N\exp(-4 \log N)= N^{-3}.
\end{equation}
In particular:
$$\mathbb{P} \big(\big| |R_N|-m_N\big|>6 \sqrt{m_N\log N}\big)= \mathbb{P}\big(|P_{N}(0)|> 6 \sqrt{m_N\log N}\big)\ll N^{-3}.
$$
 Borel-Cantelli's  lemma and the relation $\sqrt{m_N\log N}=o(m_N)$ now imply the existence of
  $\Omega_0\subset\Omega$, $\P(\Omega_0)=1$, and for every $\omega\in\Omega_0$, of an integer $N(\omega)$  such that 
both inequalities below hold:
 $${1\over2}m_N\le|R_N| \le2 m_N,\ {\rm and}\ \Vert P_N\Vert_{\infty}\ll  \sqrt{m_N\log N}, \ N\ge N(\omega).$$ 
We deduce, for $\omega\in \Omega_0$, $N\ge N(\omega)$ and any fixed $x\in \T, x\neq 0$:
$${1\over |R_N|}\,\Big|\sum_{n\le N}\delta_ne(nx)\Big|\ll C(x)/m_N$$
$$\sup_{\tau}\Big({1\over |R_N|}\,\Big|\sum_{n\le N}(\xi_n-\delta_n)e(n\tau)\Big|\Big)\ll \sqrt{\log N/m_N},$$
whence the theorem, since $m_N\to \infty$ and $\log N/m_N\to 0$. 
\end{proof}
\begin{cor} The random Furstenberg set $T$ is a.s. Hartman uniformly distributed.
\end{cor}

\noi{\bf Remarks} 1. The corollary is failing for the Furstenberg set $S$ itself although $\sigma_N\thickapprox(\log N)^2$! Indeed, we saw in Theorem \ref{fan} that, for uncountably many $x\in \T$, one has
$${1\over|S_N|}\sum_{n\in S_N}e(nx)\not\to0.$$
\noindent 2. A question naturally arises: is  Bourgain's condition (\ref{ma}) necessary and sufficient?

\subsection{Random integers and the Khinchin property}
In this subsection, we will only assume  that  condition (\ref{aa}) holds for the sequence $(\delta_k)_{k\geq 1}$, and will not need Bourgain's condition (\ref{ma}). 
We begin by a complement on Orlicz functions (defined in Section 3).

 \begin{lem}\label{montype}Let $\varphi$ be an Orlicz function such that $\psi:=\varphi\circ u$ is concave, where $u(x)=\sqrt{x}$. Let $\varepsilon=(\varepsilon_k)_{1\leq k\leq n}$ be a Rademacher sequence on a probability space $\Omega$.  If $a_1,\ldots, a_n$ are scalars, then
$$\mathbb{E}\big[\varphi(|\sum_{k=1}^n \varepsilon_k a_k|)\big]\leq \sum_{k=1}^n \varphi(|a_k|).$$
As  a corollary, if $X_1,\ldots, X_n$ are  independent, centered, $L^\varphi$-integrable functions, we have:
$$\mathbb{E}\big(\varphi(| \sum_{k=1}^n X_k|)\big)\leq  \sum_{k=1}^n \mathbb{E}\big(\varphi(2| X_k|)\big) .$$
\end{lem}
\begin{proof} Let $X=|\sum_{k=1}^n \varepsilon_k a_k|$. Jensen's inequality and the subaddivity of $\psi$ give 
$$\mathbb{E}\big[\varphi(|X|)\big]=\mathbb{E}\big[\psi(|X|^2)\big]\leq \psi\big(\mathbb{E}(|X|^2)\big)= \psi\big(\sum_{k=1}^n |a_k|^2)$$$$\leq \sum_{k=1}^n \psi(|a_k|^2)=\sum_{k=1}^n \varphi(|a_k|).$$
 And the corollary follows by symmetrization, noting that for a centered variable $Y\in L^\varphi$, it holds by Proposition \ref{raga}: 
 $$\mathbb{E}(\varphi(|Y|))\leq \mathbb{E}(\varphi(|\tilde Y|))\leq \mathbb{E}(\varphi(|2Y|)).$$
 \end{proof}
\noi {\bf Examples.} a) With $\varphi(x)=x^p,\ 1\leq p\leq 2$, we recover the fact that  the Banach space $L^p$ is of type $p$ (see e.g. \cite{LQ},  vol.1, Ch.V, p.~188).\\
b) We can also take  $\varphi(x)=\varphi_{\alpha}(x)=x\,\log^{\alpha} (1+x)=:xv(x),\ 0<\alpha\leq 1$. Let us prove the concavity of $\psi$ on $\R^+$. Since $-\frac{u''}{u'^{2}}=\frac{1}{u}$, we must show that $\varphi''(x)/\varphi'(x)\leq 1/x$, or else that 
$x^2v''(x)+xv'(x)-v(x)\leq 0$.  This is clear since $v'' \leq 0$ and $v'/v=\frac{\alpha}{(1+x)\log (1+x)}\leq\frac{1}{(1+x)\log (1+x)}\leq\frac{1}{x} $.\\

We now recall that the {\em Khinchin class} of an increasing sequence $(u_n)_{n\geq 1}$ of positive integers is the class of $L^1$-functions $f$ (here with mean zero) satisfying 
     \begin{equation}\label{tau}\tau_{k}(x):=\sum_{j=1}^k f(u_j x)=o(k) \hbox{\ almost everywhere}.\end{equation} 
Marstrand (\cite{Ma}) proved that there are  $L^\infty$-functions which are not in the Khinchin class of the set $\mathbb{N}$ of positive integers. 
Here is an interesting related result. Koksma (\cite{Ko}) proved that if the function $f=\sum_{k\neq 0} a_k e_k\in L^2$ with 
$$\sum_{|k|\geq 3} |a_k|^2 (\log\log |k|)^3<\infty,$$
that is, $\widehat{f}$ slightly better than $\ell^2$,
then $f$ belongs to the  Khinchin class of $\mathbb{N}$.
  We will prove here the following.
\begin{theo}\label{jenemar} 
Let $f$ be in the Khinchin class of the positive integers, with moreover $f\in L^{\varphi_\beta}$ for some $\beta>1$. Then, $f$ also belongs almost surely to the  Khinchin class of the random set $R=\{k \ ;\  \xi_{k}=1\}=:\{u_p\}.$ 
\end{theo}
\begin{proof}Let $E\subset \mathbb{T}$ with $m(E)=1$ and $\tau_{k}(x)=o(k)$ for $x\in E$ (see (\ref{tau})). We begin by observing that 
$$\frac{f(u_1x)+\cdots + f(u_p x)}{p} =\frac{1}{S_N}\sum_{k=1}^N \xi_k f(kx) \hbox{\ with}\ N=u_p,$$
where $S_N=\sum_{k=1}^N \xi_k=|R_N|$.
The number $N$ is random, but tends to $\infty$ almost surely, which is enough here. As we saw, almost surely, $m_N/2\leq S_N\leq 2m_N$ for $N$ large; hence, we can as well consider the quotient 
$$A_N:=\frac{1}{m_N}\sum_{k=1}^N \xi_k f(kx)=\frac{1}{m_N}{\sum_{k=1}^N (\xi_k-\delta_k) f(kx)}+\frac{1}{m_N}\sum_{k=1}^N \delta_k f(kx) $$$$=:B_N+C_N.$$
We will show that, given $x\in E$, $A_N\to 0$ almost surely. The Fubini theorem will then give the result. 
We first observe that the deterministic term $C_N$ tends to $0$, because for any sequence $(a_n)$, 
$\sum_{k=1}^n a_k =o(n)$ implies $\sum_{k=1}^n \delta_ka_k =o(\sum_{k=1}^n \delta_k)$, which is a consequence of (\ref{aa}) and an Abel summation.

\smallskip
For the random term $B_N$, we use a martingale argument.   
 Let us consider the product space 
$(\Omega \times \mathbb{T}, \mathcal{A}\times \mathcal{B}(\mathbb{T}), \P\times m)$, on which  is defined the martingale
$$
     M_n(\omega, x) =\sum_{k=1}^n \frac{(\xi_{k}(\omega) -\delta_k)f(kx)}{m_k}
$$ 
relative to the filtration $ \{\mathcal{F}_n\}$ defined by $\mathcal{F}_n=\sigma(\xi_1, \cdots, \xi_n)\times \mathcal{B}(\mathbb{T})$. Indeed
$$
    \mathbb{E}_{\omega, x} [(\xi_n -\delta_n)f(nx))|\mathcal{F}_{n-1}] = f(nx) \mathbb{E}_{\omega, x} [(\xi_n -\delta_n))|\mathcal{F}_{n-1}] =0.
$$
Now, pick $0<\alpha<\min (\beta-1, 1)$ and put $\varphi=\varphi_\alpha$. We claim that $M_n$ is bounded in $L^{\varphi}$, a fortiori in $L^1$. 
Then Doob's convergence theorem  implies the ``almost sure almost everywhere'' convergence of the series 
$$
  \sum_{k=1}^\infty \frac{(\xi_{k}(\omega) -\delta_k)f(kx)}{m_k}\cdot
$$
Applying Kronecker's lemma, we will  conclude that $B_N \to 0$.

To prove our claim, we first note the easy estimation
\begin{equation}\label{plusim} \int_{\T} |f|\log^{\alpha}\big(1+\frac{|f|}{A}\big)dm(x)\ll \frac{1}{(\log A)^{\beta-\alpha}}\hbox{\ as}\ A\to \infty.
\end{equation}
Let us fix $a\ge1$ and 
let $X_{k}(\omega)=\frac{(\xi_{k}(\omega)-\delta_k)f(kx)}{m_k}$. Note  that
$$\varphi(| X_k|/a) \leq  \frac{|\xi_k-\delta_k|}{m_k a}\, |f(kx)| \Big[\log^{\alpha} (1+ \frac{|f(kx)|}{A})\Big], 
\quad {\rm with} \ A = \frac{m_k }{|\xi_k -\delta_k|}.$$
We integrate with respect to  $x$, using (\ref{plusim}) and the invariance of $m$ under $ x\mapsto kx$, to get 
$$\int \varphi(| X_k|/a) dx \le  C \frac{|\xi_k-\delta_k|}{m_k a}\, \frac{1}{\log^{\beta-\alpha} A} =C \frac{|\xi_k-\delta_k|}{m_k a}\, \frac{1}{\log^{\beta-\alpha} \frac{m_k}{|\xi_k -\delta_k|}}.$$
We now take the expectation (w.r.~to $\omega$) and obtain
\begin{equation}\label{asur2}\mathbb{E}_{\omega, x}\big(\varphi( |X_k|/a)\big)\ll \frac{\delta_k}{m_{k}a (\log m_k)^{\beta-\alpha}}.
\end{equation}    
Noting that  the assumptions of Lemma \ref{montype}
are satisfied for fixed $x$, we get by permuting the order of integrations and changing in (\ref{asur2}) $a$ into $a/2$, $a\ge2$, 
$$\int\!\!\!\int \varphi\Big(\frac{|M_{n}(\omega, x)|}{a}\Big) d\omega dx\leq \frac{C'}{a}  \sum_{k=1}^\infty  \frac{\delta_k}{m_k(\log m_k)^{\beta-\alpha}} =: \frac{C''}{a}<\infty.$$
Indeed, for any sequence  $(a_n)$ of positive numbers, we have 
$$
\sum_{n=1}^\infty \frac{a_n}{(a_1+\cdots+a_n)\log^{\rho}(a_1+\cdots +a_n)}<\infty
$$
when $\rho>1$.
By definition of the  Orlicz norm, we get the claim: $\Vert M_n\Vert_\varphi\leq C''$.
\end{proof}
 
\subsection{$T$ is almost surely $p$-Rider}
In \cite{LQR} and \cite{LQR2}, the authors constructed inside (a variant of) $T$ a random subset \textit{of positive relative density}, which is $4/3$-Rider and $\Lambda(q)$ for each $q>2$ (We will come back to this in a forthcoming work). 
 Here, with the idea of randomly mimicking the Furstenberg set $S$, we have  to deal with the whole of $T$. We can prove the following.

\begin{theo}\label{version} Almost surely, the set $T:=T(\omega)$  is  $p$-Rider for   $p>4/3$ and not $p$-Rider for all $p<4/3$.
\end{theo}
\begin{proof}  As we already showed in Proposition \ref{mc}, the mesh condition implies that   $T$ is not  $p$-Rider for $p<4/3$. The proof of the other claim  is more elaborate, and relies on the following result,  due to L.~Rodr\'iguez-Piazza, which extends Pisier's  result for Sidon sets (\cite{Ro1}, Teorema 2.3 p.~85-86). We refer to section 3.1. for the notions involved.

\begin{theo}[\cite{Ro1}] \label{betis} Let $E$ be a set of \textnormal{positive} integers, and  $1\leq p<2$.  The following are equivalent\textnormal{:}\\     \indent {\rm 1)}  $E$ is $p$-Rider. \\
\indent {\rm 2)}   For every finite subset $A$ of $E$, there exists a subset $B$ of $A$ such that 
$$ B \hbox{\ is  quasi-independent and}\ |B|\geq \delta |A|^\varepsilon$$
 where $\varepsilon=\frac{2}{p}-1<1$  and where  $\delta$ is a positive constant.
\end{theo}
\medskip
 
We now need the following lemma (\cite{LQR}, Lemma II.1 p.115).

\begin{lem}[\cite{LQR}] \label{trio} Let $s\geq 2$ and $A$ be positive integers.  Set
$$\Omega_{s}(A)= \{\omega : T(\omega)\cap [A,\infty[ \hbox{\ contains at least a  relation of length}\ s\}.$$
 Then
 $$\mathbb{P}(\Omega_{s}(A))\leq \frac{C^s}{s^s} \sum_{j>A} \delta_{j}^2 m_{j}^{s-2}, \hbox{\ with}\ C=4e.$$
\end{lem}
We denote $T_k=T\cap [1,k],\  k=1,2,\ldots $. Then the following holds.
\begin{lem}\label{bluff} Let $A_n=\big(n\log^{3}n\big)^n$; then\\
\indent {\rm 1)} $\sum_{n\geq 1} \mathbb{P}(\Omega_{n}(A_n))<\infty.$\\
\indent {\rm 2)} Almost surely,  $|T_{A_n}|\ll (n\log n)^2$ for all integers $n $ large enough. 
\end{lem}

Lemma \ref{bluff} means that with high probability,  $T$ does not contain any ``short''  relation at infinity, and contains relatively few elements near the origin. 
\medskip

\noi {\em Proof of  Theorem \ref{version}.}
Let us take Lemma \ref{bluff}  for granted.  Denote 
$$\Omega':=\Omega\setminus\limsup_{n} \Omega_{n}(A_n),\quad 
\Omega'':=\{\omega :|T_{A_n}|\ll (n\log n)^2 \hbox{\ for}\ n\hbox{\ large}\}. $$ 
It holds $\mathbb{P}(\Omega')=\mathbb{P}(\Omega'')=1$, $\mathbb{P}(\Omega'\cap \Omega'')=1$, by the Borel-Cantelli lemma. 
Let us fix $\omega\in \Omega'\cap \Omega''$ and then $T:=T(\omega)$. Let $n_0(\omega)=:n_0\geq 1$ such that $\omega\notin \Omega_{n}(A_n)$ for $n\geq n_0$.  Then, by definition,
\begin{equation}\label{presque} \hbox{\ for}\  n\geq n_0,\ T\cap [A_n,\infty[ \hbox{\ contains no relation of length}\  \leq n.
\end{equation}
Moreover,
  \begin{equation} \label{suiv} |T_{A_n}|\ll (n\log n)^2 \hbox{\ for}\ n\geq n_1. \end{equation} 
 Now, we fix some finite subset $E$ of $T$ and we put $n=|E|^\varepsilon$. We can assume that $n\geq \max(n_0, n_1)$.  We then observe that
  $$\big|E\cap [A_n,\infty[\big|\geq \big|E\big|-\big|E\cap [1, A_n]\big|\geq n^{1/\varepsilon}-C(n\log n)^2\geq \delta n^{1/\varepsilon}\geq n$$
 for large $n$, provided that
  $1/\varepsilon>2, \hbox{\ or else}\ \varepsilon=\frac{2}{p}-1<1/2.$
  And this condition is equivalent to $p>4/3$, which was our assumption. But then, if we take  $F\subset E\cap [A_n,\infty[$ with $|F|=n$, this set is quasi-independent by (\ref {presque}), and has the right cardinality in order to apply Theorem \ref{betis}.
  
  All this shows that, for  $\omega\in \Omega'\cap \Omega''$,  $T(\omega)$ is $p$-Rider for all $p>4/3$.  This ends the proof of Theorem \ref{version} \textit{up to that of  Lemma \ref{bluff}}. 
  \end{proof}
  \noi \textit{A sketchy proof of Lemma \ref{bluff}.} The proof is sketchy since it is similar to that of Lemma II.3 in  \cite{LQR}.  First, Tchebycheff's inequality and (\ref{thha})  give  
  $$
  \mathbb{P}\Big(\big| |T_{A_n}|-\mathbb{E}(|T_{A_n}|)\big|>m_{A_n}\Big)\leq \frac{m_{A_n}}{m_{A_n}^2}=\frac{1}{m_{A_n}}\ll \frac{1}{(n\log n)^2}.
  $$
So that
 $\sum_{n\geq 1} \mathbb{P}\Big(\big| |T_{A_n}|-\mathbb{E}(|T_{A_n}|)\big|>(\log A_n)^2\Big)< \infty.$
  In other terms, one has (\ref{suiv}),  giving the second assertion.  Next, according to Lemma \ref{trio}, we have
  $$\mathbb{P}\big(\Omega_{n}(A_n)\big)\ll \frac{C^n}{n^n} \sum_{j>A_n} \frac{\log^{2}j}{j^2} (\log^{2} j)^{n-2}\ll \frac{C^n}{n^n}\int_{A_n}^\infty \frac{(\log^{2} t)^{n}}{t^2}dt$$
  $$\ll\frac{C^n}{n^n} \frac{(\log^{2} A_n)^n}{A_n}$$
  by using an integration by  parts. Or else,
  \begin{equation*}\label{adit}\mathbb{P}\big(\Omega_{n}(A_n)\big)\ll  \frac{C^n}{n^n}\frac{(n\log n)^{2n}}{(n\log^{3}n)^n}\ll \frac{C^n}{(\log n)^{n}},\end{equation*}
 which gives the first assertion.
$\hfill {\square}$


\section{$S$ and the Bohr topology}
\subsection{Reminders}
Recall that if  $G$ is a locally compact abelian group with dual $\Gamma$, the  Bohr  compactification $\beta\Gamma$ of  $\Gamma$ is the dual group of $G_d$, the group $G$ equipped with the discrete topology. The group $\beta\Gamma$ is the set of all characters (continuous or not) on $G$, it is compact and contains $\Gamma$ as a dense subgroup. One can describe more concretely  this topology when $G=\T$ and  $\Gamma=\Z$.
\begin{defi} The Bohr topology on  $\Z$ is  the  group topology with the following basis of neighbourhoods of  zero: 
$$V(x_1,\dots,x_k,\eps)=\{n\in\Z; \ |e(nx_j)-1|<\eps\,  \hbox{\ for}\ 1\le j\le k\}$$
called {\bf  Bohr neighbourhoods} (of $0$ in $\Z$).
\end{defi}
\begin{prop} 1. If $V$  is a Bohr neighbourhood of $0$, there exists another neighbourhood $W$ of $0$ such that $V\supset W-W$.\\
2. A Bohr neighbourhood  has positive upper  density.
\end{prop}
\noindent Assertion 2. follows from the pigeonhole principle and the simultaneous approximation.
\begin{prop} The Bohr topology on  $\Z$ is the coarsest  topology for which all Fourier transforms of discrete measures are continuous.
\end{prop}
Warning: this topology is non-metrizable. For example, a sequence of  distinct integers $(n_j)$ never converges in  $\beta\Z$; otherwise, $e_{n_j}(x)\to\ell(x)\not=0$ for all $x$ and 
$0=\int e_{n_{j+1}-n_j}(x) dx\to 1$ by Lebesque's convergence theorem.
We are concerned by  $\Z$ equipped with its  Bohr topology and to some subsets of integers,  dense or  not (in $\beta\Z$).

\subsection{$S$ is Bohr-closed}
%
\begin{prop} $E:=S(p_1,\dots, p_r)$, in particular $S$, is Bohr-closed. 
\end{prop}
\begin{proof} We write  $p^\alpha$ for $p_{1}^{\alpha_1}\cdots p_{r}^{\alpha_r}$, and we write  $\beta\geq \alpha$ if $\beta_j\geq \alpha_j$ for all  $j$; and  $\beta>\alpha$ if $\beta_j > \alpha_j$ for all $j$. We now show that $\Z\backslash E$ is open for the  Bohr topology.
Indeed, let $m\notin E$. We separate two cases.  \\
\noindent $\bullet$  $m=0$. Then $V=N\Z$, where $N$ has a prime factor  $>p_r$, is a neighbourhood  of  $m$ disjoint from $E$.\\
\noindent $\bullet$  $m\neq 0$. 
One  writes 
$$m=p^\alpha n=:s_0\times n$$ 
with $s_0\in E$,\  $n\not=0,1$ and $n\wedge  p_{1} p_{2}. ..p_r=1$. 
 One can find $s=p^\beta\in E$ with $\beta>0$ and the $\beta_j$'s large enough so as to have $s>n-1$, implying that   $n\not\equiv1$ mod $s$. \\
We now  claim that the neighbourhood  of $m$,
$$V:=V(m)=m+sp^\alpha \Z,$$ 
satisfies $V\cap E=\emptyset$. Indeed, a  relation
$$p^\gamma=m+sp^\alpha k= p^\alpha(n+sk),$$
clearly implies $\gamma\geq \alpha$. If $\gamma=\alpha$, then $n+sk=1$ and $n\equiv1$ mod $s$, contradicting the choice of  $s$. Therefore, we have for example  $\gamma_1>\alpha_1$. 
 After simplification by $p_{1}^{\alpha_1}$, we get
$p_1|(n+sk)$ which is again  impossible: since $p_1|s$, we have 
 $(n+sk)\wedge p_1=n\wedge p_1=1$. 
\end{proof}
\noindent {\bf Comment.} In particular, the set $S$ is not  an intersective set,  which answers negatively to a question of Bergelson.\\
\noindent We already mentioned that a $H$-ud set is Bohr-dense. Using this fact, we were able to prove that such Bohr-dense sets,  moreover $p$-Sidon for all $p>1$, do exist \cite{LQR}.  
This echoes an old, still open, conjecture:\\
``Is there a Sidon(=$1$-Sidon) set $E\subset \mathbb{Z}$ which is Bohr-dense?''
The question is still widely open, and a new approach will be needed, because\\ 
$\bullet$ A Sidon set $E$ is never H-ud.\\
This can be seen for example by Hartman's version of  Drury's theorem for Sidon sets \cite{HAR}: write $E=(u_n)_{n\geq 1}$; there exists a \textit{continuous} measure $\mu$ on $\T$ such that
$$\widehat{\mu}(u_n)=-1 \hbox{\ for all}\ n\in \N.$$ 
But then $\sigma=\delta_0+\mu$ satisfies $\widehat{\sigma}(u_n)=0 \hbox{\ for}\ n\in \N$; so that
$$A_N:=\int_{\T}\frac{1}{N} \sum_{n=1}^N e(-u_{n}x)d\sigma(x)=\frac{1}{N} \sum_{n=1}^N \widehat{\sigma}(u_n)=0.$$
But if  $E$ is $H$-ud, $A_N\to \sigma(\{0\})=1$ by Lebesgue's theorem, and we have a contradiction.\\
  Another approach would be as follows: if $E$ is $H$-ud, then $C_E$ contains a copy of the Banach space $c_0$ of sequences which tend to zero at infinity \cite{PIQ}. But since $E$ is Sidon, $C_E\sim \ell_1$. And notably, $\ell_1$ does not contain $c_0$. \\
$\bullet$ For a Sidon set $E$, $W(E)$ is uncountable. More generally, it has been observed by Hartman that $W(E)$ is uncountable as soon as some continuous measure exists with $\inf_{n\in E}|\widehat\mu(n)|\ge\delta>0$. Indeed, if the assertion is false with $E=(n_k)$, we then have
$$\lim_N{1\over N}\sum_{n\le N} e(n_kt)\to 0$$ for every $t$ outside some countable set. For any continuous measure $\nu$, we get by integration
$$\lim_N{1\over N}\sum_{n\le N} \widehat\nu(n_k)\to 0;$$
now, by taking $\nu=\mu\ast\tilde\mu$ where $\mu$ satisfies $|\widehat\mu(n)|>\delta$ on $E$, we get a contradiction.
\medskip

\subsection*{Acknowledgments.} The authors warmly thank the referee for a very careful reading, which allowed  a substantially improved presentation of the paper.
A. H. Fan is partially supported by NSFC (grant no.11971192). 
 H.~Queff\'elec and M.~Queff\'elec acknowledge support of the Labex CEMPI (ANR-11-LABX-0007-01).


\end{document}